\pdfoutput=1
\RequirePackage{ifpdf}
\ifpdf 
\documentclass[pdftex]{sigma}
\else
\documentclass{sigma}
\fi

\usepackage{tikz-cd}

\numberwithin{equation}{section}

\newtheorem{Theorem}{Theorem}[section]
\newtheorem*{Theorem*}{Theorem}
\newtheorem{Corollary}[Theorem]{Corollary}
\newtheorem{Lemma}[Theorem]{Lemma}
\newtheorem{Proposition}[Theorem]{Proposition}
 { \theoremstyle{definition}
\newtheorem{Definition}[Theorem]{Definition}

\newtheorem{Example}[Theorem]{Example}
\newtheorem{Remark}[Theorem]{Remark} }

\DeclareMathOperator{\Span}{Span}

\DeclareMathOperator{\fix}{Fix}

\DeclareMathOperator{\Imag}{Im}
\DeclareMathOperator{\Real}{Re}

\DeclareMathOperator{\Hom}{Hom}

\DeclareMathOperator{\Gr}{Gr}

\DeclareMathOperator{\Tr}{Trace}

\newcommand{\orb}{\mathcal{O}}

\newcommand{\R}{\mathbb{R}}
\newcommand{\C}{\mathbb{C}}
\newcommand{\CP}{\mathbb{C}P}
\newcommand{\D}{\mathcal{D}}
\newcommand{\CS}{\mathbb{CS}}
\newcommand{\RP}{\mathbb{R}P}
\newcommand{\Q}{\mathbb{H}}

\newcommand{\Id}{{\rm Id}}

\begin{document}
\allowdisplaybreaks

\newcommand{\arXivNumber}{2009.10417}

\renewcommand{\PaperNumber}{114}

\FirstPageHeading

\ShortArticleName{Real Forms of Holomorphic Hamiltonian Systems}

\ArticleName{Real Forms of Holomorphic Hamiltonian Systems}

\Author{Philip ARATHOON~$^{\rm a}$ and Marine FONTAINE~$^{\rm b}$}

\AuthorNameForHeading{P.~Arathoon and M.~Fontaine}

\Address{$^{\rm a)}$~Department of Mathematics, University of Michigan, Ann Arbor, MI 48109, USA}
\EmailD{\href{mailto:philash@umich.edu}{philash@umich.edu}}

\Address{$^{\rm b)}$~Mathematics Institute, University of Warwick, Coventry, CV4 7AL, UK}
\EmailD{\href{mailto:marine.fontaine@warwick.ac.uk}{marine.fontaine@warwick.ac.uk}}

\ArticleDates{Received June 12, 2024, in final form December 10, 2024; Published online December 21, 2024}

\Abstract{By complexifying a Hamiltonian system, one obtains dynamics on a holomorphic symplectic manifold. To invert this construction, we present a theory of real forms which not only recovers the original system but also yields different real Hamiltonian systems which share the same complexification. This provides a notion of real forms for holomorphic Hamiltonian systems analogous to that of real forms for complex Lie algebras. Our main result is that the complexification of any analytic mechanical system on a Grassmannian admits a real form on a compact symplectic manifold. This produces a `unitary trick' for Hamiltonian systems which curiously requires an essential use of hyperk\"{a}hler geometry. We demonstrate this result by finding compact real forms for the simple pendulum, the spherical pendulum, and the rigid body.}


\Classification{53D20; 14J42}

\section{Background and outline}
Real analytic Hamiltonian systems are closely related to complex holomorphic Hamiltonian systems. Indeed, if we treat the variables in a real analytic system as being complex we obtain a~holomorphic system. For this reason, it is not uncommon to treat these two concepts equivalently. However, what this perspective overlooks is the possibility that many different and distinct real analytic systems might each complexify to the same system. This motivates a theory of real forms for holomorphic Hamiltonian systems which will allow us to treat different real Hamiltonian systems as real forms of the same complex system.

This idea is not new and has been considered before in \cite{bulgarian3,bulgarian4,bulgarian1,bulgarian2}. However, these previous approaches are limited to systems on $\C^{2n}$ with real subspaces $\R^{2n}$ as real forms. We~generalise this to dynamics on holomorphic symplectic manifolds and introduce a wider definition for what it means to be a real form. This has the advantage of extending the scope of the theory to include a greater variety of dynamical systems and also brings it into closer contact with ideas in differential geometry. Below we give an outline of the work contained within.

\subsection{Introduction}

A real analytic manifold can be complexified to give a complex manifold \cite{kulkarni,realform2,realform1}. The original manifold appears as a totally real submanifold of half dimension, and thus, any such submanifold with this property shall be considered a real form of the complex manifold. If~the~original manifold possesses an analytic symplectic form, then the complexification will be a holomorphic symplectic manifold. Exactly as with ordinary Hamiltonian dynamics, one can consider a Hamiltonian vector field generated by a holomorphic function and investigate the dynamics generated by its flow.

We distinguish those real forms upon which the holomorphic symplectic form is either purely real or pure imaginary; what we call a real- or imaginary-symplectic form. If such a real form is invariant under the flow of a holomorphic Hamiltonian, then we justify why the restricted dynamics on this real form can be considered to be a real form of the holomorphic Hamiltonian system.

\subsection{Reduction}
If a holomorphic Hamiltonian systems admits some symmetry, it is natural to ask how this symmetry might manifest on a real form. A more precise formulation of this question is to ask how the complex and real symplectic reduced spaces might be related. We choose to address this in terms of Poisson reduction. For when the symmetry of the holomorphic symplectic manifold is in a certain sense compatible with respect to a real form, we show the extent to which the Poisson reduced space for the real form can itself be considered a real form of the holomorphic Poisson reduced space.

To demonstrate these results, we introduce a guiding example which finds use throughout the paper: the example of a particular coadjoint orbit ${\rm Orb}(\zeta)$ in $\mathfrak{gl}_n\C^*$. Geometrically, this orbit is the complex symmetric space $\D_\C$ of decompositions of $\C^n$ into two complementary subspaces. Using a dual pair construction, we exhibit this orbit as a reduced space for an action of ${\rm GL}_m\C$ on $T^*_{1,0}\Hom(\C^m,\C^n)$. We then find various real forms inside $T^*_{1,0}\Hom(\C^m,\C^n)$ which are compatible with respect to the group action, and hence, descend to real-symplectic forms on $\D_{\C}$. These turn out to be related to the real ${\rm GL}_n\R$, unitary ${\rm U}(l,n-l)$, and quaternionic ${\rm GL}_{n/2}\Q$ real forms of ${\rm GL}_n\C$.

\subsection{Branes}
A hyperk\"{a}hler manifold defines a holomorphic symplectic form for each choice of complex structure. Naturally, this provides us with a stock of examples of holomorphic symplectic manifolds. In addition to this, complex-Lagrangian submanifolds taken with respect to one complex structure give real- and imaginary-symplectic forms with respect to two other complex structures, respectively. This gives us a nice trick to generate examples of real- and imaginary-symplectic forms. We adopt the terminology taken from ideas in string theory and refer to such submanifolds as branes \cite{witten2, witten1}.

We implement this trick for a particular manifold which we obtain through hyperk\"{a}hler reduction for an action of ${\rm U}(m)$ on $\Hom(\Q^m,\Q^n)$. For different choices of complex structure, the underlying holomorphic symplectic manifold is either the coadjoint orbit ${\rm Orb}(\zeta)$ or the cotangent bundle to a complex Grassmannian. In this way, we demonstrate that complex-Lagrangian submanifolds of $T^*_{1,0}\Gr_\C$, such as the zero section or a fibre, correspond to real-symplectic forms of $\D_\C$. We also obtain an explicit diffeomorphism between $\D_{\C}$ and $T^*_{1,0}\Gr_\C$ which interchanges their hyperk\"{a}hler structures and turns out to be an essential ingredient for the compact real-form result in the following section.

Incidentally, we generate new examples of complex-Lagrangian submanifolds of $T^*_{1,0}\Gr_\C$ and provide symplectomorphisms between certain symmetric spaces with cotangent bundles of Grassmannians.

\subsection{Integrability and dynamical systems}

The generalisation of an integrable system to the holomorphic category is straightforward, considered for instance in \cite{pcrook} and extensively discussed for the real analytic category in \cite{Moerbeke, realm}. Our main question is to ask whether holomorphic integrability is equivalent to integrability on a real form. We answer in the affirmative and show that if a real form admits an analytic integrable system, then this may be extended to a holomorphic integrable system in a neighbourhood of the real form. Conversely, if a holomorphic Hamiltonian system is integrable in the holomorphic sense, then any real form of the dynamics is integrable in the real sense. This provides us with a method for generating new integrable systems from old: begin with a real analytic integrable system, complexify it and obtain real integrable systems on different real forms. We are particularly interested in when the real forms are compact symplectic manifolds. In this regard, we show the following.
\begin{Theorem*}
For any real analytic mechanical system defined on a real, complex, or quaternionic Grassmannian there exists a real Hamiltonian system defined on a non-empty open subset of a~compact symplectic manifold for which both systems complexify to the same holomorphic Hamiltonian system.
\end{Theorem*}
We demonstrate this result for the examples of the simple pendulum, the spherical pendulum, and the rigid body. The phases spaces for these problems are the cotangent bundles $T^*S^1$, $T^*S^2$, and $T^*{\rm SO}(3)$, and their compact real forms turn out to be $S^2$, $S^2\times S^2$, and $\CP^3$, respectively.

\section{Introduction}
\subsection{Holomorphic symplectic geometry}
A holomorphic symplectic manifold is a complex manifold $M$ equipped with a closed and non-degenerate holomorphic 2-form $\Omega$. This gives an isomorphism between the holomorphic tangent bundle $T^*_{1,0}M$ and the real tangent bundle $TM$ by identifying the complex-valued covector $df$ to $x\in M$ with the tangent vector $X_f$ satisfying
\[
\Omega(X_f,Y)=\langle \mathrm{d}f, Y\rangle
\]
for all $Y\in T_xM$. In this way, we may associate to a holomorphic function $f$ a Hamiltonian vector field $X_f$.
If we separate the holomorphic symplectic form into its real and imaginary parts we obtain real symplectic forms $\omega_R$ and $\omega_I$ on $M$, where $\Omega=\omega_R+{\rm i}\omega_I$. If we write the complex structure on each tangent space as $I$, then these two forms are related by
\begin{equation}
	\label{wRwIrelated}
	\omega_R(I(X),Y)=-\omega_I(X,Y).
\end{equation}
\begin{Proposition}\label{bihamiltonian}
	Let $f$ be a holomorphic function defined on a holomorphic symplectic manifold~$(M,\Omega)$ and write $f=u+{\rm i}v$ as its decomposition into real and imaginary parts. The~Hamiltonian vector field $X_f$ is equivalently the Hamiltonian vector field of $u$ on $(M,\omega_R)$ and of~$v$ on~$(M,\omega_I)$.
\end{Proposition}
\begin{proof}
	This follows immediately by expanding $\Omega(X_f,Y)=\langle \mathrm{d}f, Y\rangle$ into real and imaginary parts as
	\[
	\omega_R(X_f,Y)+{\rm i}\omega_I(X_f,Y)=\langle \mathrm{d}u, Y\rangle+{\rm i}\langle \mathrm{d}v,Y\rangle.
\tag*{\qed}
 \]
\renewcommand{\qed}{}
\end{proof}

This proposition shows that the Hamiltonian vector field generated by a holomorphic function satisfies a bihamiltonian structure on $M$ with respect to the two symplectic forms $\omega_R$ and $\omega_I$.

To generalise the concept of real structures on complex vector spaces, we say that a real structure on a complex manifold $M$ is an involution $R$ whose derivative is everywhere conjugate-linear. If non-empty, the fixed-point set of $R$ can be considered a real form of $M$. More generally, a real form of a complex manifold shall mean a totally real submanifold of half dimension. Recall that a submanifold $N\subset M$ is called totally real if $T_xN\cap I(T_xN)=\{0\}$ for all $x$ in $N$.

We immediately caution that our use of the term `real form' should not be confused with the concept of a real-valued differential form. Our choice of language is designed to reflect the corresponding notion of real forms in Lie theory.

\begin{Definition}\label{first_defn}
	A real form $N$ of $(M,\Omega)$ is called a \emph{real-symplectic form} if the restriction of~$\Omega$ to~$N$ is purely real, and an \emph{imaginary-symplectic form} if this restriction is purely imaginary. Furthermore, a real structure $R$ is said to be
	\begin{itemize}\itemsep=0pt
		\item a \emph{real-symplectic structure} if $R^*\Omega=\overline{\Omega}$ (conjugate-symplectic),
		\item and an \emph{imaginary-symplectic structure} if $R^*\Omega=-\overline{\Omega}$ (anti-conjugate-symplectic).
	\end{itemize}
\end{Definition}
Assuming they are non-empty, the fixed-point sets of real- and imaginary-symplectic structures are real- and imaginary-symplectic forms, respectively.
\begin{Proposition}\label{symp_Lang}
	A totally real submanifold $N$ of $(M,\Omega)$ is a real-symplectic form if and only if it is a Lagrangian submanifold with respect to $(M,\omega_I)$. This implies $N$ is a symplectic sub\-manifold of $(M,\omega_R)$. Likewise, $N$ is an imaginary-symplectic form if and only if it is Lagrangian with respect to $(M,\omega_R)$, and this implies it is a symplectic submanifold of $(M,\omega_I)$.
\end{Proposition}
\begin{proof}
	The first statement is immediate from the definitions. If $T_xN$ is a Lagrangian subspace with respect to $\omega_I$, then from \eqref{wRwIrelated} it follows that $I(T_xN)$ is the orthogonal complement to $T_xN$ with respect to $\omega_R$. As $I(T_xN)\cap T_xN=\{0\}$, we see that $T_xN$ is a symplectic subspace of $T_xM$ with respect to $\omega_R$ for all $x\in N$. The case for imaginary-symplectic forms is similar.
\end{proof}

\begin{Example}[cotangent lift of a real structure]
	Let $C$ be a complex manifold together with a real structure $r$ with fixed-point set $\fix r=C^r$. We can lift this to a real structure $R_\pm$ on the holomorphic cotangent bundle by setting
	\[
	\langle R_\pm(\eta),X\rangle=\pm\overline{\langle \eta,r_*X\rangle}
	\]
	for $\eta\in T^*_{1,0}C$ a covector to $x$ and for all $X\in T_{r(x)}C$. This satisfies $R_\pm^*\lambda=\pm\overline{\lambda}$, where $\lambda$ is the canonical one-form. Therefore, we may lift real structures on $C$ to either real- or imaginary-symplectic structures on $(T^*_{1,0}C,\Omega_{can})$. The fixed-point set $\fix R_+$ is canonically symplectomorphic to $T^*C^r$, and $\fix R_-$ to ${\rm i}T^*C^r$, by which we mean the bundle of imaginary-valued 1-forms on $C^r$.
\end{Example}

Suppose we have a dynamical system on $(M,\Omega)$ generated by a holomorphic Hamiltonian $f$. Given a real form $N\subset M$ we would like to be able to describe a real Hamiltonian system on $N$ which can, in a sense, be said to be a real form of the holomorphic system $(M,\Omega,f)$. This raises two questions: what real forms $N$ should we consider, and what should be the corresponding real Hamiltonian?

In answer to the first question, we shall insist that $N$ is a real- or imaginary-symplectic form. For the sake of brevity, we shall throughout this paper mostly make reference to real-symplectic forms, however the imaginary case is entirely similar. From Proposition~\ref{symp_Lang}, the restriction of~$\omega_R$ to $N$ is symplectic. If we write this restriction as~$\widehat{\omega}_R$, then we have the real symplectic manifold~$(N,\widehat{\omega}_R)$.

We now turn to the second question. In order for the dynamical system on $M$ to induce a~dynamical system on $N$, we must suppose that $N$ is invariant under the flow generated by~$f$. If~we~decompose the function into its real and imaginary parts as $f=u+{\rm i}v$, then by Proposition~\ref{bihamiltonian} the flow of $f$ is equivalently the Hamiltonian flow of $u$ on $(M,\omega_R)$. For $x\in N$, the Hamiltonian vector field $X_u$ belongs to $T_xN$ if and only if $du$ yields zero when evaluated on the orthogonal complement to $T_xN$ with respect to $\omega_R$. From the proof of Proposition~\ref{symp_Lang}, this complement is $I(T_xN)$, and therefore we require $\langle \mathrm{d}u, I(T_xN)\rangle=0$. Since $f$ is holomorphic, the Cauchy--Riemann equations show us that this is equivalent to $\langle \mathrm{d}v, T_xN\rangle=0$.
\begin{Theorem}Let $(M,\Omega)$ be a holomorphic symplectic manifold. A real-symplectic form ${N\subset M}$ is invariant under the Hamiltonian flow generated by a holomorphic function $f=u+{\rm i}v$ if and only if $v$ is locally constant on $N$. In this case, the flow on $N$ is identical to the Hamiltonian flow generated by the restriction of $u$ to $(N,\widehat{\omega}_R)$.
\end{Theorem}
For the situation described in this theorem, we can speak of the dynamical system $(N,\widehat{\omega}_R,u)$ as being a real form of the holomorphic system $(M,\Omega,f)$. The principal application we have in mind is to identify different real Hamiltonian systems by recognising that they are both real forms of the same holomorphic system.

More generally, we are interested in holomorphic Hamiltonians whose restriction to $N$ is real up to a constant phase. After all, we may always scale such a Hamiltonian by the constant phase to obtain a real-valued restriction which generates a Hamiltonian system on $N$. For this reason, in the presence of a real-/imaginary-structure we have the following.

\begin{Definition}\label{compat_Ham}
	Let $(M,\Omega)$ be a holomorphic symplectic manifold and $R$ a real-/imaginary-symplectic structure. A holomorphic Hamiltonian $H$ will be called $R$-compatible with respect to a real structure $*$ on $\C$ if
	\begin{equation*}
		H\circ R=H^*.
	\end{equation*}
\end{Definition}

\subsection{Holomorphic Poisson geometry}

A holomorphic Poisson manifold is a complex manifold $P$ equipped with a holomorphic section~$\pi$ of $\bigwedge^{2,0}TP$ with the property that the bracket
\[
\{f,g\}(x)=\pi_x(\mathrm{d}f,\mathrm{d}g)
\]
defined on the sheaf of holomorphic functions satisfies the Jacobi identity and is thus a com\-plex-valued Poisson bracket. Such a Poisson structure defines a fibrewise complex-linear map $\sharp\colon T^*_{1,0}P\rightarrow TP$ which sends a covector $\mathrm{d}f\in T^*_{1,0}P$ at $x$ to the tangent vector $\sharp \mathrm{d}f\in T_xP$ satisfying
\begin{equation}
	\label{sharpdefn}
	\langle \mathrm{d}g,\sharp \mathrm{d}f\rangle=\pi_x(\mathrm{d}f,\mathrm{d}g)
\end{equation}
for all $\mathrm{d}g\in T_{1,0}^*P$ at $x$. Exactly as with the real situation, the image of $\sharp$ defines a complex involutive (generalised) distribution called the characteristic distribution. One may define a~holomorphic symplectic form $\Omega$ on each leaf $\mathcal{O}$ of the distribution by
\begin{equation}\label{induced_symp_form}
	\Omega(\sharp \mathrm{d}f,\sharp \mathrm{d}g)=\pi(\mathrm{d}f,\mathrm{d}g).
\end{equation}
This form is non-degenerate and well defined thanks to $\ker\sharp=(\Imag\sharp)^\circ$ and is a closed holomorphic form as a consequence of the Jacobi identity. The leaves of the characteristic distribution are therefore immersed holomorphic symplectic manifolds.

By decomposing a complex one-form into its real and imaginary parts $\mathrm{d}f=\mathrm{d}u+{\rm i}\mathrm{d}v$, we establish two real isomorphisms $\mathrm{d}f\leftrightarrow \mathrm{d}u$ and $\mathrm{d}f\leftrightarrow \mathrm{d}v$ between the space of complex-linear forms $T^*_{1,0}P$ with real-linear forms $T^*P$. We can then define two real Poisson structures $\pi_R$ and~$\pi_I$ on $P$ by setting
\begin{equation*}
	\pi(\mathrm{d}f_1,\mathrm{d}f_2)=\pi_R(\mathrm{d}u_1,\mathrm{d}u_2)+{\rm i}\pi_I(\mathrm{d}v_1,\mathrm{d}v_2).
\end{equation*}
Applying \eqref{sharpdefn} gives
\begin{equation}\label{sharprelations}
	\sharp \mathrm{d}f=\sharp_R \mathrm{d}u=\sharp_I \mathrm{d}v.
\end{equation}
Consequently, the characteristic distributions coincide, and with the aid of~\eqref{induced_symp_form} we see that the real and imaginary parts of the holomorphic symplectic form $\Omega$ on a leaf $\mathcal{O}$ are precisely the real symplectic forms $\omega_R$ and $\omega_I$ induced by the Poisson structures $\pi_R$ and $\pi_I$, respectively.

In light of Proposition~\ref{symp_Lang}, the generalisation of Definition~\ref{first_defn} requires a review of the appropriate analogues of Lagrangian and symplectic submanifolds in symplectic geometry to Poisson geometry.
\begin{itemize}\itemsep=0pt
	\item A submanifold $N$ of $(P,\pi_I)$ is coisotropic if $\sharp_I(TN^\circ)\subset TN$.
	
	\item A submanifold $N$ of $(P,\pi_R)$ is called pointwise Poisson--Dirac if $\sharp_R(T_xN^\circ)\cap T_xN=\{0\}$ for every $x\in N$. For any $\mathrm{d}U\in T_x^*N$, this implies the existence of a unique $\mathrm{d}u\in T_x^*P$ which projects to $\mathrm{d}U$ and for which $\sharp_R(\mathrm{d}u)\in T_xN$. In this way, we may define a bivector~$\Pi_R$ on~$N$ by
	\begin{equation}
		\label{Bivector}
		\Pi_R(\mathrm{d}U_1,\mathrm{d}U_2)=\pi_R(\mathrm{d}u_1,\mathrm{d}u_2).
	\end{equation}
	If $\Pi_R$ varies smoothly across $N$, then $N$ is a Poisson--Dirac submanifold of $(P,\pi_R)$.
\end{itemize}

\begin{Definition}
	A real form $N\subset P$ is a \emph{real-Poisson form} of $(P,\pi)$ if $N$ is a coisotropic submanifold of $(P,\pi_I)$ and an \emph{imaginary-Poisson form} if $N$ is a coisotropic submanifold of~$(P,\pi_R)$. A real structure $R$ on a holomorphic Poisson manifold will be called a \emph{real-} ($+$) or \emph{imaginary-}~($-$) \emph{Poisson structure} if
	\begin{equation}\label{real_Poisson_structure}
		\pi_{R(x)}\bigl(\overline{R^*\mathrm{d}f},\overline{R^*\mathrm{d}g}\bigr)=\pm~\overline{\pi_{x}(\mathrm{d}f,\mathrm{d}g)}
	\end{equation}
	holds at all $x$ and for all holomorphic one-forms $\mathrm{d}f,\mathrm{d}g\in T_{1,0}^*P$ at $x$. The notation $\overline{R^*\mathrm{d}f}$ denotes the conjugate-adjoint, which for all $X\in T_{R(x)}P$ satisfies
	\begin{equation}
		\label{conjugate_adjoint}
		\bigl\langle \overline{R^*\mathrm{d}f},X\bigr\rangle=\overline{\langle \mathrm{d}f,R_*X\rangle}.
	\end{equation}
	
\end{Definition}
If non-empty, the fixed-point set $N$ of a real-Poisson structure $R$ is a real-Poisson form. To~see this, consider the complex-linear form $\mathrm{d}f=\mathrm{d}u+{\rm i}\mathrm{d}v$ for $\mathrm{d}v\in T_xN^\circ$. For $X$ in $T_xN$, we have that $\langle \mathrm{d}f, X\rangle$ is real and $\langle \mathrm{d}f,I(X)\rangle$ is imaginary, from which we see that $\overline{R^*\mathrm{d}f}=\mathrm{d}f$ at $x$ in $N$. It follows from \eqref{real_Poisson_structure} that $\pi(\mathrm{d}f_1,\mathrm{d}f_2)$ is purely real, and so $0=\pi_I(\mathrm{d}v_1,\mathrm{d}v_2)=\langle \mathrm{d}v_2,\sharp_I \mathrm{d}v_1\rangle$ for all $\mathrm{d}v_1, \mathrm{d}v_2\in T_xN^\circ$. The vector $\sharp_I \mathrm{d}v_1$ must therefore belong to $T_xN$, and hence, $N$ is a coisotropic submanifold of $(P,\pi_I)$.

\begin{Proposition}If $N$ is a real-Poisson form of a holomorphic Poisson manifold $(P,\pi)$, then~$N$ is a Poisson--Dirac submanifold of $(P,\pi_R)$. Furthermore, if the intersection between $N$ and a~holomorphic symplectic leaf $\orb$ is a submanifold of $\orb$, then $N\cap\orb$ is a real-symplectic form of~$(\orb,\Omega)$.
\end{Proposition}
\begin{proof}
If $N$ is a real-Poisson form, then the subspace $T_xN\cap T_x\orb$ is a coisotropic subspace of~$(T_x\orb,\omega_I)$. On the other hand, since $N$ is a real form and $\orb$ a complex submanifold, we must have
	\[
	(T_xN\cap T_x\orb)\oplus I(T_xN\cap T_x\orb)\subset T_x\orb.
	\]
	Therefore, the dimension of $T_xN\cap T_x\orb$ must be less than or equal to half the dimension of~$T_x\orb$. Yet since this is a coisotropic subspace it must be exactly half the dimension, and hence, a~Lagrangian subspace. It follows from Proposition~\ref{symp_Lang} that if $N\cap\orb$ is a submanifold of $\orb$, then it is Lagrangian with respect to $(\orb,\omega_I)$ and therefore a real-symplectic form of $(\orb,\Omega)$.
	
	Since $T_xN\cap T_x\orb$ is Lagrangian with respect to $(T_x\orb,\omega_I)$, it follows from \eqref{wRwIrelated} that it is symplectic with respect to $(T_x\orb,\omega_R)$ and therefore, $N$ is pointwise Poisson--Dirac with respect to~$(P,\omega_R)$. To show that $\Pi_R$ in \eqref{Bivector} varies smoothly, it suffices to show that the map $T^*N\rightarrow T^*P|_N$ which sends $\mathrm{d}U\in T_x^*N$ to $\mathrm{d}u\in T_x^*P$ is smooth \cite{crainic}. As $T_xP=T_xN\oplus I(T_xN)$, we may smoothly extend $\mathrm{d}U$ to a form on $T_xP$ by setting it to equal zero on $I(T_xN)$. We claim that this extension is precisely $\mathrm{d}u$. By complex linearity of the form $\mathrm{d}f=\mathrm{d}u+{\rm i}\mathrm{d}v$, we see that $\mathrm{d}u$ vanishing on $I(T_xN)$ implies $\mathrm{d}v$ vanishes on $T_xN$, and hence $\mathrm{d}v\in T_xN^\circ$. Let $\mathrm{d}w\in T_xN^\circ$ be arbitrary. From~\eqref{sharprelations}, we obtain $\langle \mathrm{d}w, \sharp_R \mathrm{d}u\rangle=\langle \mathrm{d}w, \sharp_I \mathrm{d}v\rangle$ which must equal zero as $\sharp_I \mathrm{d}v\in T_xN$ since $N$ is coisotropic, and so $\sharp_R \mathrm{d}u\in T_xN$ as desired.
\end{proof}

\begin{Example}[complex Lie algebras as holomorphic Poisson manifolds and their real forms]\label{proto}
	The prototypical example of a holomorphic Poisson manifold is the dual of a complex Lie algebra~$\mathfrak{g}^*$ equipped with the Kostant--Kirilov--Souriau (KKS) Poisson bracket
	\[
	\{f,g\}(\eta)=\pi_\eta(\mathrm{d}f|_\eta,\mathrm{d}g|_\eta)=\langle \eta, [\mathrm{d}f|_\eta,\mathrm{d}g|_\eta]\rangle,
	\]
	where the one-forms $\mathrm{d}f|_\eta$ and $\mathrm{d}g|_\eta$ on $\mathfrak{g}^*$ belong to the Lie algebra~$\mathfrak{g}$ upon which the Lie bracket~$[\,,\,]$ is defined. Let $\rho_*$ be a real form on $\mathfrak{g}$ which is also a Lie algebra automorphism with non-empty fixed-point set $\mathfrak{g}^\rho$. In other words, $\mathfrak{g}^\rho$ is a real form of $\mathfrak{g}$ in the Lie algebraic sense. This involution lifts to a real structure $\overline{\rho^*}$ on $\mathfrak{g}^*$ given by the conjugate-adjoint as in \eqref{conjugate_adjoint}. This defines a real-Poisson structure on $\mathfrak{g}^*$ whose fixed-point set may be identified with the dual of~$\mathfrak{g}^\rho$. The negative $-\overline{\rho^*}$ defines an imaginary-Poisson structure whose fixed-point set is the space of imaginary-valued 1-forms on $\mathfrak{g}^\rho$.
\end{Example}

\section{Reduction}
\subsection{Holomorphic Poisson reduction}
Let $(M,\Omega)$ be a holomorphic symplectic manifold and $G$ a complex Lie group which acts on $M$ by holomorphic symplectomorphisms. The quotient topology on the orbit space $M/G$ is not nice in general, so we shall suppose that the orbit map $P\colon M\rightarrow M/G$ is a holomorphic submersion between complex manifolds. By virtue of $G$ acting symplectically, the Poisson bracket between $G$-invariant functions is again $G$-invariant. This allows us to define a unique Poisson structure $\widetilde{\pi}$ on $M/G$ for which the projection map is a Poisson map. The space $(M/G,\widetilde{\pi})$ is the (holomorphic) {Poisson reduced space}.

Consider the fixed-point set $M^R$ of a real-symplectic structure $R$ on $M$. We would like to understand how this behaves with respect to Poisson reduction. In the presence of a real structure it is reasonable to expect some degree of compatibility between the involution $R$ and the action of $G$. If we suppose that $R$ maps $G$-orbits into $G$-orbits, then it descends to a real structure on $M/G$ which we shall denote by $\widetilde{R}$. Using $\widetilde{R}\circ P=P\circ R$, it follows from \eqref{real_Poisson_structure} and the definition of $\widetilde{\pi}$ that $\widetilde{R}$ is a real-Poisson structure on $(M/G,\widetilde{\pi})$. This compatibility condition can be ensured if the following equivariant definition holds.
\begin{Definition}
	A holomorphic group action of a complex Lie group $G$ on a complex manifold~$M$ equipped with a real structure $R$ will be called \emph{$R$-compatible} with respect to a real group structure $\rho$ on $G$ if for all $g\in G$ and $x\in M$
	\begin{equation}
		\label{equivariant_condition}
		R(g\cdot x)=\rho(g)\cdot R(x).
	\end{equation}
\end{Definition}

For such a compatible group action, the real Lie group $G^\rho$ acts symplectically on the real form $\bigl(M^R,\widehat{\omega}_R\bigr)$. If in addition, we suppose the orbit map $p\colon M^R\rightarrow M^R/G^\rho$ is also a submersion between smooth manifolds, then we can equally consider the real Poisson reduced space \smash{$\bigl(M^R/G^\rho,\widehat{\Pi}_R\bigr)$}. The following proposition establishes the relation between the two possible choices of real reduced space, \smash{$\bigl(M^R/G^\rho,\widehat{\Pi}_R\bigr)$} and \smash{$\bigl((M/G)^{\widetilde{R}},\widetilde{\Pi}_R\bigr)$}.

\begin{Theorem}\label{reductionthm}
	There is a Poisson map $\Psi$ from $\bigl(M^R/G^\rho,\widehat{\Pi}_R\bigr)$ into $\bigl((M/G)^{\widetilde{R}},\widetilde{\Pi}_R\bigr)$. This map is an immersion with discrete fibres and image $P\bigl(M^R\bigr)$. Moreover, if $G$-acts freely on $M^R$, then this map is an injection, and hence, $P\bigl(M^R\bigr)$ is an immersed Poisson submanifold of \smash{$(M/G)^{\widetilde{R}}$}.
\end{Theorem}

\begin{proof}\samepage
	The map $\Psi$ sends the $G^\rho$-orbit through $x\in M^R$ to the $G$-orbit through $x$. The commutativity of the square below tells us that the image of $\Psi$ is $P\bigl(M^R\bigr)$.
	\begin{equation*}
			\begin{tikzcd}
				M^R\arrow[r, hook,"\iota"]\arrow[d,two heads,"p"'] & M\arrow[d,two heads, "P"] \\ M^R/G^\rho \arrow[r,hook,"\Psi"] & M/G
			\end{tikzcd}
			\end{equation*}
	
Since $P\circ\iota=\Psi\circ p$ is smooth, it follows from an application of the submersion theorem that~$\Psi$ is a smooth map.
	
For $x$ in $M^R$, consider the tangent vector $\xi\cdot x$ for $\xi\in\mathfrak{g}$. If $\xi\cdot x$ belongs to $T_xM^R$, then $R_*(\xi\cdot x)=\xi\cdot x$. By taking the infinitesimal version of \eqref{equivariant_condition}, we have $\xi\cdot x=\rho_*(\xi)\cdot x$. It follows that $\xi\cdot x$ is also the tangent vector generated by the element $(\xi+\rho_*(\xi))/2$. However, this element is clearly fixed by $\rho_*$ and so belongs to the Lie algebra $\mathfrak{g}^\rho$ of $G^\rho$. It follows that the tangent vector $\xi\cdot x$ belongs to $T_x(G^\rho\cdot x)$, which establishes
	\begin{equation}\label{likeOshea}
		T_xM^R\cap T_x(G\cdot x)=T_x(G^\rho\cdot x).
	\end{equation}
	Consequently, the intersection of a $G$-orbit with $M^R$ is the discrete union of $G^\rho$-orbits and so the fibres of $\Psi$ are discrete. For $x\in M^R$, suppose $g\cdot x$ also belongs to $M^R$ for some $g\in G$. By~\eqref{equivariant_condition}, we have $g\cdot x=\rho(g)\cdot x$, and therefore, if $G$ acts freely on $M^R$, then $\rho(g)=g$. Therefore, $(G\cdot x)\cap M^R=G^\rho\cdot x$, which implies that $\Psi$ is an injection.
	
	Let $x(t)$ be a curve in $M^R$ with tangent vector $X$ at $x=x(0)$ and suppose $\Psi_*p_*X$ is zero in $T_{P(x)}(M/G)$. This implies $X\in T_x(G\cdot x)\cap T_xM^R$ which from \eqref{likeOshea} shows that $p_*X=0$, and hence, $\Psi$ is an immersion.
	
	Let $U_1$ and $U_2$ be locally defined real functions on $(M/G)^{\widetilde{R}}$. Since this is a Poisson--Dirac submanifold of $(M/G,\widetilde{\pi}_R)$ there exist extensions $u_1$ and $u_2$ on $M/G$ whose Hamiltonian vector fields on $(M/G)^{\widetilde{R}}$ are tangent to $(M/G)^{\widetilde{R}}$ and for which
	\begin{equation}\label{interm}
		\{U_1,U_2\}_{\widetilde{\Pi}_R}(P(x))=\{u_1,u_2\}_{\widetilde{\pi}_R}(P(x))=\{u_1\circ P,u_2\circ P\}_{\pi_R}(x)
	\end{equation}
	for any $x\in M^R$. In the last equality we have used the fact that $P$ is a Poisson map. The extensions $u_1$ and $u_2$ may be assumed to be $\widetilde{R}$-invariant. As $P\circ R=\widetilde{R}\circ P$, the functions $u_1\circ P$ and $u_2\circ P$ are also $R$-invariant, and therefore their Hamiltonian vector fields are tangent to $M^R$. Since $M^R$ is a symplectic submanifold with respect to $\omega_R$, the right-hand side above is equal to
	\[
	\{u_1\circ(P\circ\iota),u_2\circ(P\circ\iota)\}_{\widehat{\pi}_R}(x),
	\]
	where $\widehat{\pi}_R$ is the Poisson bivector for $\bigl(M^R,\widehat{\omega}_R\bigr)$. As $P\circ\iota=\Psi\circ p$, we can use the fact that $p$ is Poisson to rewrite this as
	\[
	\{u_1\circ\Psi,u_2\circ\Psi\}_{\widehat{\Pi}_R}(p(x)).
	\]
	By comparing this to the left-hand side in \eqref{interm} and writing $P(x)=\Psi(p(x))$, it follows from surjectivity of $p$ that $\Psi$ is Poisson.
\end{proof}

\begin{Remark}If the action of $G$ is free and admits an equivariant holomorphic momentum map, then the holomorphic symplectic leaves of $(M/G,\widetilde{\pi})$ are connected components of the orbit-reduced spaces for the Hamiltonian action of $G$ on $M$. If in addition, the $G$-action is $R$-compatible with respect to a real form $\rho$, then the action of $G^\rho$ on $M^R$ is also free and Hamiltonian and the symplectic leaves of $\bigl(M^R/G^\rho,\widehat{\Pi}_R\bigr)$ are the connected components of the orbit-reduced spaces. The previous theorem tells us that $\Psi(M^R/G^\rho)$ is an immersed Poisson submanifold of \smash{$\bigl((M/G)^{\widetilde{R}},\widetilde{\Pi}_R\bigr)$}. This implies that $\Psi\bigl(M^R/G^\rho\bigr)$ is a union of symplectic leaves in \smash{$(M/G)^{\widetilde{R}}$}, and since $\Psi$ is Poisson, it restricts to a symplectomorphism between leaves in $M^R/G^\rho$ and leaves in \smash{$(M/G)^{\widetilde{R}}$}. It follows that $\Psi$ restricted to a $G^\rho$-orbit-reduced space on $M^R$ gives a~symplectomorphism between this space and a real-symplectic form of a $G$-orbit-reduced space~on~$M$.
\end{Remark}

\subsection{A coadjoint orbit example}\label{theexample}
Consider the space $\Hom(\C^m,\C^n)$ of complex $n\times m$ matrices $Q$. The dual space is identified with matrices $P$ in $\Hom(\C^n,\C^m)$ via the trace pairing $\Tr(PQ)$. The holomorphic cotangent bundle $T^*_{1,0}\Hom(\C^m,\C^n)$ is then identified with the set of such pairs $(Q,P)$. The holomorphic symplectic form is
\begin{equation*}
	\Omega((Q_1,P_1),(Q_2,P_2))=\Tr(P_2Q_1-P_1Q_2).
\end{equation*}
The groups ${\rm GL}_m\C$ and ${\rm GL}_n\C$ act symplectically by $\bigl(Qg^{-1},gP\bigr)$ and $\bigl(hQ,Ph^{-1}\bigr)$ with momentum maps
\begin{equation*}	
	\mu_m(Q,P)=PQ,\qquad\text{and}\qquad\mu_n(Q,P)=QP.
\end{equation*}
Here we have identified $\mathfrak{gl}_n$ and $\mathfrak{gl}_m$ with their respective duals using the trace form.

Now suppose $m\le n$ and let $M$ be the open subset of the cotangent bundle where both $Q$ and $P$ have maximal rank. One can show (see \cite{dualpairs}) that $\mu_n$ is an orbit-map for the ${\rm GL}_m\C$-action on $M$ and that the fibres of $\mu_m$ are orbits of ${\rm GL}_n\C$. As the momentum map is Poisson we can identify the Poisson reduced space $M/{\rm GL}_m\C$ with $\mu_n(M)\subset\mathfrak{gl}_n\C^*$. For $\zeta\in\C^\times$, the orbit-reduced space $\mu_m^{-1}(\zeta\cdot\Id_m)/{\rm GL}_m\C$ is identified with the coadjoint orbit
\begin{equation}
	\label{coad_orbit}
	{\rm Orb}(\zeta)=\{\xi=QP\mid (Q,P)\in M,\,PQ=\zeta\cdot\Id_m\}\subset\mathfrak{gl}_n\C^*.
\end{equation}

Geometrically, $\xi$ determines a decomposition $\C^n=\Sigma^{(m)}\oplus \Lambda^{(n-m)}$, where ${\Sigma=\Imag Q}$ and $\Lambda=\ker P$. Conversely, such a decomposition determines an element $\xi$ where $\ker\xi=\Lambda$ and ${\xi|_\Sigma=\zeta\cdot\Id_\Sigma}$. Naturally, we see that the orbit is diffeomorphic to the symmetric space $\D_\C$ given~by
\begin{Definition}
	For $\mathbb{F}=\R,\C,\Q$, the \emph{space of decompositions} $\D_\mathbb{F}$ is the set
	\begin{equation}
		\label{elegant_A+B}
		\D_\mathbb{F}=(\Gr_\mathbb{F}(m;n)\times \Gr_\mathbb{F}(n-m;n))\setminus\Delta,
	\end{equation}
	where $\Delta$ is the closed subset of pairs $(\Sigma,\Lambda)$ with $\Sigma\cap\Lambda\ne\{0\}$.
\end{Definition}
\begin{Remark}\label{elegant_rmk}
	The symplectic form on $\D_{\C}$ admits a rather elegant description. Recall that tangent vectors to an element $\Sigma$ of the Grassmannian may be identified with linear maps ${\Sigma\rightarrow\C^n/\Sigma}$. Since $(\Sigma,\Lambda)$ determines a decomposition of $\C^n$, the tangent space to $\D_{\C}$ at this point may be identified with pairs of maps $(\alpha\colon\Sigma\rightarrow\Lambda,\,\beta\colon \Lambda\rightarrow\Sigma)$. The KKS form on ${\rm Orb}(\zeta)$ can be shown~to~be
	\[
		\Omega_{\rm KKS}((\alpha_1,\beta_1),(\alpha_2,\beta_2))=\zeta\Tr(\beta_2\alpha_1-\beta_1\alpha_2).
	\]
\end{Remark}

\subsection{Real forms of the coadjoint orbit}

We shall now exhibit real-symplectic forms of ${\rm Orb}(\zeta)$ via two simultaneous methods: as fixed-point sets arising from real group structures on ${\rm GL}_n\C$ as in Example~\ref{proto}, and using real-symplectic structures on $M$ which descend through the reduction as in Theorem~\ref{reductionthm}.

The group ${\rm GL}_n\C$ has three types of real forms up to isomorphism:
\begin{enumerate}\itemsep=0pt
	\item Real: Equip $\C^n$ with a real structure $x\mapsto\overline{x}$. The real linear group ${\rm GL}_n\R$ is the fixed-point set of $\rho(g)=\overline{g}$.
	\item Unitary: Equip $\C^n$ with a non-degenerate Hermitian form of signature $(l,n-l)$ given by $(x,y)=x^\dagger \Theta y$. The unitary group ${\rm U}(l,n-l)$ is the fixed-point set of $\sigma(g)=\Theta^{-1}g^{-\dagger}\Theta$.
	\item Quaternionic: For $n$ even, endow $\C^n$ with a quaternionic structure by equipping it with a conjugate-linear map $\mathbb{J}$ with $\mathbb{J}^2=-\Id$. The quaternionic linear group ${\rm GL}_{n/2}\Q$ is the fixed-point set of $\tau(g)=-\mathbb{J}g\mathbb{J}$.
\end{enumerate}
The final ingredient we need is an additional quaternionic structure on $\C^m$ in the case $m$ is even which we also denote by $\mathbb{J}$. With all of these structures in place we can introduce the following involutions on $T^*_{1,0}\Hom(\C^m,\C^n)$:
\begin{equation}\label{list_of_RST}
	\begin{aligned}
		&R_\pm(Q,P)=(\overline{Q},\pm\overline{P}), \\
		&T_\pm(Q,P)=(-\mathbb{J}Q\mathbb{J},\mp\mathbb{J}P\mathbb{J}),
	\end{aligned}
	\qquad
	\begin{aligned}
		&S_+(Q,P)=\bigl({\rm i}\Theta^{-1}P^\dagger, iQ^\dagger \Theta\bigr),\\
		&S_-(Q,P)=\bigl(-\Theta^{-1}P^\dagger, -Q^\dagger \Theta\bigr).
	\end{aligned}
\end{equation}
These are real-symplectic structures for $+$ and imaginary-symplectic structures for $-$. The ${\rm GL}_m\C$-action is $R_\pm$-, $S_\pm$-, and $T_\pm$-compatible with respect to $\rho$, $\sigma$, and $\tau$, respectively. Moreover, they each descend under the quotient to give the real-Poisson structures $\pm\overline{\rho^*}$, $\pm\overline{\sigma^*}$, and $\pm\overline{\tau^*}$ appearing in Example~\ref{proto}.

\begin{Definition}The symmetric space $\D_{\C}^{\rm unit}\subset\D_{\C}$ is the set of decompositions $\C^n=\Sigma\oplus\Lambda$ which are orthogonal with respect to the Hermitian form $(x,y)=x^\dagger \Theta y$.
\end{Definition}
If the Hermitian form is positive-definite, $\D_{\C}^{\rm unit}\cong\Gr_\C(m;n)$. On the other hand, if the form is indefinite, $\D_{\C}^{\rm unit}$ contains multiple connected components each corresponding to the signature of the Hermitian form restricted to the component spaces of the decomposition.

\begin{Proposition}\label{RST}
	For $\zeta$ real, $\overline{\rho^*}$ and $\overline{\tau^*}$ restrict to real-symplectic structures on ${\rm Orb}(\zeta)$, and $-\overline{\sigma^*}$ to an imaginary-symplectic structure. Conversely, for $\zeta$ imaginary, $-\overline{\rho^*}$ and~$-\overline{\tau^*}$ restrict to imaginary-symplectic structures, and $\overline{\sigma^*}$ to a real-symplectic structure. By identifying ${\rm Orb}(\zeta)$ with $\D_{\C}$, the fixed-point sets of these structures are
	\begin{equation*}
		\fix(\pm\overline{\rho^*})=\D_\R,\qquad\fix(\pm\overline{\sigma^*})=\D_{\C}^{\rm unit},\qquad\fix(\pm\overline{\tau^*})=\D_\Q.
	\end{equation*}
\end{Proposition}
\begin{proof}
	This can be shown by verifying that the involutions in \eqref{list_of_RST} descend through the orbit map $\mu_n$ to give the involutions $\overline{\rho^*}(\xi)=\overline{\xi}$, $\overline{\sigma^*}(\xi)=-\Theta^{-1}\xi^\dagger\Theta$, and $\overline{\tau^*}(\xi)=-\mathbb{J}\xi\mathbb{J}$.
\end{proof}

\begin{Remark}
	Let $\langle x,y\rangle=x^\dagger y$ denote the ordinary Hermitian forms on $\C^n$ and $\C^m$. This allows us to equip $T^*_{1,0}\Hom(\C^m,\C^n)$ with a K\"{a}hler metric
	\begin{equation}\label{metric}
		g((Q_1,P_1),(Q_2,P_2))=\frac{1}{2}\Tr\bigl(Q_1Q_2^\dagger+Q_2Q_1^\dagger\bigr)+\frac{1}{2}\Tr\bigl(P_1P_2^\dagger+P_2P_1^\dagger\bigr).
	\end{equation}
	Looking ahead, we would like the involutions in \eqref{list_of_RST} to be isometries of this metric. To ensure this, we will from now impose the following:
	\begin{itemize}\itemsep=0pt
		\item Compatibility between the real structure on $\C^n$ with the standard Hermitian form $\langle \,,\,\rangle$, the~Hermitian form $(\,,\,)$, and the quaternionic structures $\mathbb{J}$, in the sense that: $\langle\overline{x},\overline{y}\rangle=\overline{\langle x, y\rangle}$, $(\overline{x},\overline{y})=\overline{(x,y)}$, and $\mathbb{J}(\overline{x})=\overline{\mathbb{J}(x)}$.
		\item Compatibility between the Hermitian form $(\,,\,)$ defined by $\Theta$ and quaternionic structures~$\mathbb{J}$ with $\langle\,,\,\rangle$ in the sense that: $\langle \Theta x,\Theta y\rangle=\langle x,y\rangle$ and $\langle\mathbb{J}(x),\mathbb{J}(y)\rangle=\langle x,y\rangle$.
	\end{itemize}
	In plainer terms, if $x\mapsto\overline{x}$ is ordinary complex conjugation on $\C^n$, we require $\Theta$ to be a real symmetric matrix with $\Theta^2=\Id_n$, and $\mathbb{J}(x)=\mathcal{J}\overline{x}$ for $\mathcal{J}$ a real skew-symmetric matrix.
\end{Remark}

\section{Branes}
\subsection{Hyperk\"{a}hler geometry}
A hyperk\"{a}hler manifold is a Riemannian manifold $(M,g)$ equipped with three complex structures, $I$, $J$, and $K$, which satisfy the usual quaternionic relations, and for which $(M,g)$ is K\"{a}hler with respect to each of them. A hyperk\"{a}hler manifold defines a holomorphic symplectic form on $M$ with respect to each complex structure. If we denote the K\"{a}hler forms by
\[
\omega_1(X,Y)=g(I(X),Y),\qquad\omega_2(X,Y)=g(J(X),Y),\qquad\omega_3(X,Y)=g(K(X),Y),
\]
then $\Omega_1=\omega_2+{\rm i}\omega_3$, $\Omega_2=\omega_3+{\rm i}\omega_1$, and $\Omega_3=\omega_1+{\rm i}\omega_2$ each define holomorphic symplectic forms on $M$ with respect to the complex structures, $I$, $J$, and $K$, respectively.

\begin{Proposition}\label{branes}
	A submanifold $N$ of a hyperk\"{a}hler manifold $M$ is a complex-Lagrangian submanifold of $(M,I,\Omega_1)$ if and only if it is an imaginary-symplectic form of $(M,J,\Omega_2)$ and a~real-symplectic form of $(M,K,\Omega_3)$.
\end{Proposition}
\begin{proof}
	It suffices to consider the tangent space $T_xM$ to a point $x\in N$. If $N$ is complex-Lagrangian with respect to $(M,I,\Omega_1)$, then it is Lagrangian with respect to $\omega_2$ and $\omega_3$. This~implies that $g(J(X),Y)$ and $g(K(X),Y)$ are zero for all $X,Y\in T_xN$, which means $J(T_xN)=K(T_xN)$ is the orthogonal complement to $T_xN$. Consequently, $T_xN$ is a real form with respect to both $J$ and $K$, and so the first implication follows from Proposition~\ref{symp_Lang}.
	
	Conversely, if $N$ is an imaginary-symplectic form of $(M,J,\Omega_2)$ and a real-symplectic form of $(M,K,\Omega_3)$, then it is Lagrangian with respect to both $\omega_3$ and $\omega_2$. This implies $J(T_xN)$ and $K(T_xN)$ are both equal to the orthogonal complement to $T_xN$, and so $N$ is a complex submanifold of $(M,I)$ since $I(T_xN)=JK(T_xN)=T_xN$.
\end{proof}

\begin{Definition}
	A \emph{brane} of a hyperk\"{a}hler manifold shall mean a submanifold which is complex-Lagrangian with respect to some holomorphic symplectic form associated to the hyperk\"{a}hler structure.
\end{Definition}

\subsection{A hyperk\"{a}hler structure on the coadjoint orbit}
Consider the quaternionic spaces $\Q^m$ and $\Q^n$ viewed as right $\Q$-modules. Then the space $\Hom(\Q^m,\Q^n)$ of $n\times m$ quaternionic matrices is itself a left $\Q$-module, where we denote left multiplication by the quaternions ${\rm i}$, ${\rm j}$, ${\rm k}$
with the operators $I$, $J$, $K$. We can further equip this space with a metric%
\begin{equation}\label{quat_metric}
	g(A,B)=\frac{1}{2}\Tr\bigl(AB^\dagger+BA^\dagger\bigr),
\end{equation}
where $A^\dagger$ is the quaternionic conjugate-transpose. The tuple $(g,I,J,K)$ defines a hyperk\"{a}hler structure on $\Hom(\Q^m,\Q^n)$.
Observe that right multiplication by the compact symplectic group ${\rm Sp}(m)\subset {\rm GL}_m\Q$ preserves the hyperk\"{a}hler structure. The following lemma will facilitate us in switching between different complex structures on this space.
\begin{Lemma}\label{facilitate}
	By writing elements $A$ of $\Hom(\Q^m,\Q^n)$ as
	\begin{equation*}
		A=\begin{cases}
			\displaystyle Q+{\rm j}P^\dagger,\\[1mm]
			\displaystyle\frac{Y^\dagger-X}{\sqrt{2}}-{\rm j}\biggl(\frac{{\rm i}\overline{X}+{\rm i}Y^{\mathsf{T}}}{\sqrt{2}}\biggr),\vspace{1mm}\\
			\displaystyle\frac{U+{\rm i}V^\dagger}{\sqrt{2}}+{\rm j}\biggl(\frac{{\rm i}V^\dagger-U}{\sqrt{2}}\biggr)
		\end{cases}
	\end{equation*}
	for $(Q,P)$, $(X,Y)$, and $(U,V)$ pairs of matrices in $\Hom(\C^m,\C^n)\times\Hom(\C^n,\C^m)$, we establish isometric holomorphic symplectomorphisms between $T^*_{1,0}\Hom(\C^m,\C^n)$ and $\Hom(\Q^m,\Q^n)$ for the structures $(I,\Omega_1)$, $(J,\Omega_2)$, and $(K,\Omega_3)$, respectively. Furthermore, these intertwine the ${\rm U}(m)\subset {\rm Sp}(m)$-action with the ${\rm U}(m)\subset {\rm GL}_m\C$-action on $T^*_{1,0}\Hom(\C^m,\C^n)$.
\end{Lemma}
\begin{proof}
	It is straightforward to verify that the metric in \eqref{metric} is pulled back to \eqref{quat_metric}. One can then show that for the three choices of variables above, the roles of $I$, $J$, and $K$ are cyclically permuted, and thus, so are the holomorphic symplectic forms $\Omega_1$, $\Omega_2$, and $\Omega_3$.
\end{proof}

From now on, suppose $m\le n$ and let $M$ denote the open subset of $\Hom(\Q^m,\Q^n)$ consisting of all elements $A=Q+{\rm j}P^\dagger$ for which $Q$ has maximal rank. The ${\rm U}(m)$-action on $M$ admits a~triple of momentum maps $\mu_k\colon M\rightarrow\mathfrak{u}(m)^*$ for each K\"{a}hler form $\omega_k$. This allows us to introduce the hyperk\"{a}hler quotient
\[
\tilde{M}=\bigl(\mu_1^{-1}({\rm i}\cdot\Id_m)\cap\mu_2^{-1}(0)\cap\mu_3^{-1}(0)\bigr)/{\rm U}(m).
\]
\begin{Theorem}\label{HK_red_thm}
	Let $\bigl(\tilde{M},g,I,J,K\bigr)$ be the hyperk\"{a}hler reduced space for the action of ${\rm U}(m)\subset {\rm Sp}(m)$ on $M\subset\Hom(\Q^m,\Q^n)$ at the regular value $({\rm i}\cdot\Id_m,0,0)\in\mathfrak{u}(m)^*\otimes\R^3$. We have the~following holomorphic symplectomorphisms:
	\begin{itemize}\itemsep=0pt
		\item $\bigl(\tilde{M}, I, \Omega_1\bigr)\cong (T^*_{1,0}\Gr_\C, {\rm i}, \Omega_{\rm can})$,
		\item $\bigl(\tilde{M},J,\Omega_2\bigr)\cong({\rm Orb}(-1), {\rm i}, \Omega_{\rm KKS})$,
		\item $\bigl(\tilde{M},K,\Omega_3\bigr)\cong({\rm Orb}({\rm i}), {\rm i}, \Omega_{\rm KKS})$.
	\end{itemize}
	Here ${\rm Orb}(\zeta)$ is the coadjoint orbit in $\mathfrak{gl}_n\C^*$ equipped with the Kostant--Kirilov--Souriau form $\Omega_{\rm KKS}$ and $\Gr_\C$ is the Grassmannian of complex $m$-dimensional subspaces in $\C^n$.
\end{Theorem}

\begin{proof}
	For the distinguished complex structure $I$, we use Lemma~\ref{facilitate} to identify $\Hom(\Q^m,\Q^n)$ with ${T^*_{1,0}\Hom(\C^m,\C^n)}$. The ${\rm U}(m)$-action admits a holomorphic extension to the ${\rm GL}_m\C$-action demonstrated in Section~\ref{theexample} with momentum map $\mu_m=\mu_2+{\rm i}\mu_3$ given by $\mu_m(Q,P)=PQ$.
	
	In this basis $\mu_1(Q,P)={\rm i}\bigl(Q^\dagger Q-PP^\dagger\bigr)$. We claim that each ${\rm GL}_m\C$-orbit in $M$ intersects $\mu_1^{-1}({\rm i}\cdot\Id_m)$ precisely in a single orbit of ${\rm U}(m)$. To show this, consider how the action of $g$ in ${\rm GL}_m\C$ sends $Q^\dagger Q-PP^\dagger$ to
	\begin{equation}\label{svd}
		g^{-\dagger}Q^\dagger Qg^{-1}-gPP^\dagger g^{\dagger}.
	\end{equation}
	If $Q$ has maximal rank, we may act by the appropriate $g$ to assume that $Q^\dagger Q$ is the identity. If~we~then use the singular value decomposition $g=uDv^\dagger$ and suppose that $v^\dagger PP^\dagger v$ is a~diagonal matrix $\Lambda$ with non-negative diagonal entries, then \eqref{svd} becomes $u\bigl(D^{-\dagger}D^{-1}-D\Lambda D^\dagger\bigr)u^\dagger$. We~can always find a $D$ for which this equals the identity, and hence, have established that the intersection is non-empty.
	
	To show that this intersection is a single ${\rm U}(m)$-orbit, let $\mu_1(Q,P)={\rm i}\cdot\Id_m$ and suppose that~\eqref{svd} is equal to $\Id_m$. By again writing $g$ as $uDv^\dagger$, we find from substituting ${Q^\dagger Q\!=\Id_m\!+\!PP^\dagger}$ into \eqref{svd} that we must have $D^{-\dagger}(\Id_m+S)D^{-1}=\Id_m+DSD^\dagger$, where $S=v^\dagger PP^\dagger v$. This can only hold when $D$ is the identity, and hence, $g$ is unitary.
	
	From this claim, it follows that $\bigl(\tilde{M},I,\Omega_1\bigr)$ may be identified with the holomorphic symplectic reduced space
	\begin{equation*}
		\mu_m^{-1}(0)/{\rm GL}_m\C.
	\end{equation*}
	There is a well-defined identification between this quotient and $T^*_{1,0}\Gr_\C$ which sends the equivalence class $[(Q,P)]$ to the plane $\Pi=\Imag Q$ and covector $QP$, interpreted as a map $\C^n/\Pi\rightarrow\C^n$. Observe that the canonical one-form on the cotangent bundle pulls back under this identification to the canonical one-form on $T^*_{1,0}\Hom(\C^m,\C^n)$.
	
	For a different choice of distinguished complex structure, we apply the same argument. For~$J$, we must consider $\mu_m^{-1}(-\Id_m)$ for $\mu_m=\mu_3+{\rm i}\mu_1$, and for $K$ we have $\mu_m^{-1}({\rm i}\cdot\Id_m)$ where ${\mu_m=\mu_1+{\rm i}\mu_2}$. In both cases, the momentum conditions imply that the matrices $(X,Y)$ and~$(U,V)$ all have maximal rank, and a slight alteration to the claim above establishes that the quotients may be identified with the holomorphic symplectic reduced spaces in \eqref{coad_orbit}.
\end{proof}

\begin{Remark}
	The theorem establishes that $T^*_{1,0}\Gr_\C$ and ${\rm Orb}(\zeta)$ admit a hyperk\"{a}hler metric. This metric on $T^*_{1,0}\Gr_\C$ is the Calabi metric appearing in \cite{calabi_metric} whose K\"{a}hler form we denote by~$\omega_\text{C}$. The K\"{a}hler form on ${\rm Orb}(\zeta)$ we write as $\omega_\text{K}$. Following on, from Remark~\ref{elegant_rmk} one can show that on ${\rm Orb}(\zeta)\cong\D_{\C}$ the K\"{a}hler form is given by
	\begin{equation*} \omega_\text{K}((\alpha_1,\beta_1),(\alpha_2,\beta_2))=\frac{{\rm i}|\zeta|^2}{2}\Tr\bigl(\alpha_1^\dagger\alpha_2-\alpha_2^\dagger\alpha_1\bigr)+\frac{{\rm i}|\zeta|^2}{2}\Tr\bigl(\beta_1^\dagger\beta_2-\beta_2^\dagger\beta_1\bigr).
	\end{equation*}
\end{Remark}

\subsection{A curious diffeomorphism}

Theorem~\ref{HK_red_thm} allows us to identify ${\rm Orb}(\zeta)$ with $T^*_{1,0}\Gr_\C$. This identification is quite curious as it is not a biholomorphism with respect to their standard complex structures. This has been explored before in \cite{french}. For our purposes we will need to consider the geometry of this identification in greater detail.

We begin by claiming that the identification is equivariant with respect to the coadjoint action of ${\rm U}(n)\subset {\rm GL}_n\C$ on ${\rm Orb}(\zeta)\subset\mathfrak{gl}_n\C^*$ and with respect to the cotangent lift of ${\rm U}(n)$ to $T^*_{1,0}\Gr_\C$. To see this, notice that in addition to the ${\rm Sp}(m)$-action on $\Hom(\Q^m,\Q^m)$ given by multiplication to the right, there is also an action of ${\rm Sp}(n)$ given by left-multiplication. In order to preserve the left $\Q$-module structure on $\Hom(\Q^m,\Q^m)$ this action is left matrix multiplication but with the element-wise multiplication occurring to the right. Using Lemma~\ref{facilitate}, this ${\rm U}(n)\subset {\rm Sp}(n)$-action can be seen to intertwine with the ${\rm U}(n)\subset {\rm GL}_n\C$-action on $T^*_{1,0}\Hom(\C^m,\C^m)$ for any of the distinguished complex structures $I$, $J$, or $K$. By navigating our way through the proof of Theorem~\ref{HK_red_thm}, we find that this action descends through the reduction to give the desired actions on ${\rm Orb}(\zeta)$ and $T^*_{1,0}\Gr_\C$.

\begin{Theorem}
	There exists a ${\rm U}(n)$-equivariant isometric diffeomorphism
	\begin{equation*}
		\Phi\colon \ {\rm Orb}({\rm i})\longrightarrow T^*_{1,0}\Gr_\C
	\end{equation*}
	with the properties
	\begin{gather}
			\Phi^*\Real\Omega_{\rm can}=\Imag\Omega_{\rm KKS},\qquad
			\Phi^*\Imag\Omega_{\rm can}=\omega_\textup{K},\qquad
			\Phi^*\omega_\textup{C}=\Real\Omega_{\rm can}.\label{symp_forms}
	\end{gather}
	Write $\Phi(\xi)=(\Pi,\eta)$, where $\Pi$ is an element of $\Gr_\C$ and $\eta\colon\C^n/\Pi\rightarrow\Pi$ is a covector to $\Pi$. The map $\Phi$ is given explicitly by
	\begin{align}
		&\Pi=\bigl\{x+\sqrt{\xi^\dagger \xi}(x)\mid x\in\Imag\xi\bigr\}, \label{Pi_explicit}\\
		&\eta=\frac{1}{2}\bigl(\sqrt{\xi^\dagger \xi}-\sqrt{\xi\xi^\dagger}-{\rm i}\xi-{\rm i}\xi^\dagger\bigr).\label{eta_explicit}
	\end{align}
\end{Theorem}
\begin{proof}
	The map $\Phi$ is obtained by identifying $\tilde{M}$ with ${\rm Orb}({\rm i})$ and $T^*_{1,0}\Gr_\C$. We have already argued above why this identification is ${\rm U}(n)$-equivariant. The pullbacks in \eqref{symp_forms} follow by keeping track of the K\"{a}hler forms $\omega_1$, $\omega_2$, $\omega_3$ and recognising that $\Omega_\text{can}\equiv\Omega_1$ and $\Omega_{\rm KKS}\equiv\Omega_3$.
	
	To describe the map $\Phi$ explicitly, we must delve back into the proof of Theorem~\ref{HK_red_thm}. The~element $(\Pi,\eta)$ is represented by the equivalence class of ${\rm U}(m)$-orbits through $(Q,P)$ satisfying the three momentum conditions
	\[\mu_1(Q,P)={\rm i}\bigl(Q^\dagger Q-PP^\dagger\bigr)={\rm i}\cdot\Id_m\qquad\text{and}\qquad (\mu_2+{\rm i}\mu_3)(Q,P)=PQ=0.\]
	Using Lemma~\ref{facilitate}, we may write
	\begin{equation}\label{QP_UV}
		Q=\frac{U+{\rm i}V^\dagger}{\sqrt{2}}\qquad\text{and}\qquad P=-\frac{U^\dagger+{\rm i}V}{\sqrt{2}}
	\end{equation}
	to see that the corresponding element $\xi$ of ${\rm Orb}({\rm i})$ is represented by the equivalence class of ${\rm U}(m)$-orbits through $(U,V)$ satisfying the same three momentum conditions
	\begin{equation*}
		\mu_3(U,V)={\rm i}\bigl(U^\dagger U-VV^\dagger\bigr)=0\qquad\text{and}\qquad(\mu_1+{\rm i}\mu_2)(U,V)=VU={\rm i}\cdot\Id_m.
	\end{equation*}
	The subspace $\Pi$ is the image of $Q$, and hence, is equal to the set $\bigl\{U(y)+{\rm i}V^\dagger(y)\mid y\in\C^m\bigr\}$. Using the momentum condition $V U(y)={\rm i}y$ and writing $x= U(y)$, this set can be written as~$\bigl\{x+V^\dagger V(x)\mid x\in\Imag U\bigr\}$. Recall from the proof of Theorem~\ref{HK_red_thm} that $\xi=UV$, and so the image of $U$ is also the image of $\xi$. Furthermore, the $\mu_3$-momentum condition tells us that~${U^\dagger U=VV^\dagger}$, from which it follows that $\xi^\dagger \xi$ is equal $\bigl(V^\dagger V\bigr)^2$. These are all positive-semidefinite Hermitian matrices, meaning the matrix square root is well defined, from which \eqref{Pi_explicit} now follows.
	
Finally, we recall once more from the proof of Theorem~\ref{HK_red_thm} that $\eta=QP$, and then use \eqref{QP_UV} to write this in terms of $U$ and $V$. Another application of the condition $U^\dagger U=VV^\dagger$ combined with the matrix square root gives \eqref{eta_explicit}.	 	
\end{proof}

We are especially interested in the case $m=1$ for which $\Gr_\C$ is the projective space $\CP^{n-1}$. Recall from \eqref{elegant_A+B} that ${\rm Orb}({\rm i})$ may be identified with the space of decompositions of $\C^n$ into a~line $\Sigma$ and a complementary $(n-1)$-dimensional space $\Lambda$. Using the standard Hermitian form on $\C^n$, we may take the orthogonal complement $\Lambda^\perp$ and instead view the orbit as the set of pairs of non-orthogonal lines $\bigl(\Sigma,\Lambda^\perp\bigr)$.
\begin{Proposition}\label{Phi_for_CPn}
	For when $m=1$, the map $\Phi\colon {\rm Orb}({\rm i})\rightarrow T^*_{1,0}\Gr_\C$ sending $\xi$ to $(\Pi,\eta)$ may be identified with the map
	\[
	\bigl(\CP^{n-1}\times\CP^{n-1}\bigr)\setminus\Delta\longrightarrow T^*_{1,0}\CP^{n-1},
	\]
	which sends a pair of non-orthogonal lines $\Sigma=\Imag\xi=\Span\{x\}$ and $\Lambda^\perp=(\ker\xi)^\perp=\Span\{y\}$ to the line $\Pi=\Span\{z\}$ and the covector
	\begin{equation*}
		\eta=\frac{|x||y|}{2\langle y,x\rangle}zw^\dagger.
	\end{equation*}
	Here we have introduced
	\begin{equation}\label{z_w_defn}
		z=\hat{x}+\frac{\langle y, x\rangle}{|\langle x,y\rangle|}\hat{y}\qquad\text{and}\qquad w=\hat{y}-\frac{\langle x,y\rangle}{|\langle x,y\rangle|}\hat{x}.
	\end{equation}
\end{Proposition}
\begin{proof}
	For the pair $\bigl(\Sigma,\Lambda^\perp\bigr)$, the corresponding element $\xi$ of ${\rm Orb}({\rm i})$ is the map with $\ker\xi=\Lambda$ and $\xi|_\Sigma={\rm i}\cdot\Id_\Sigma$. We therefore have
	\begin{equation}\label{xi_as_xy}
		\xi={\rm i}\frac{xy^\dagger}{\langle y,x\rangle},
	\end{equation}
	from which we obtain
	\[
	\sqrt{\xi^\dagger\xi}=\frac{|x|}{|y||\langle x,y\rangle|}yy^\dagger\qquad\text{and}\qquad\sqrt{\xi\xi^\dagger}=\frac{|y|}{|x||\langle x, y\rangle|}xx^\dagger.
	\]
	The expressions for $\Pi$ and $\eta$ now follow immediately from \eqref{Pi_explicit} and \eqref{eta_explicit}.
\end{proof}

\subsection{Branes of the coadjoint orbit}
Theorem~\ref{HK_red_thm} together with Proposition~\ref{branes} allow us to identify various branes of $\tilde{M}$ with submanifolds of $T^*_{1,0}\Gr_\C$ and ${\rm Orb}(\zeta)$. The rows in the table below serve as a mnemonic to keep track of the types of these submanifolds,
\begin{equation}\label{little_table}
	\begin{tabular}{c|c|c}
		$T^*_{1,0}\Gr_\C$ & ${\rm Orb}(-1)$ & ${\rm Orb}({\rm i})$\\\hline
		\text{real} & complex & imaginary \\\hline
		complex & imaginary & real \\\hline
		imaginary & real & complex
	\end{tabular}
\end{equation}
For instance, the first row says that complex-Lagrangian submanifolds of ${\rm Orb}(-1)$ are identified with real-symplectic forms of $T^*_{1,0}\Gr_\C$ and imaginary-symplectic forms of ${\rm Orb}({\rm i})$. We already have a few examples of real- and imaginary-symplectic forms of ${\rm Orb}(\zeta)$ using Proposition~\ref{RST}. Our goal in this section is to show that these examples provide branes in $\tilde{M}$ and to describe the corresponding submanifolds which appear in the table above.

In a slight abuse of notation, to denote the involutions on $\Hom(\Q^m,\Q^n)$ we write $R_\pm$, $S_\pm$, and $T_\pm$ given by
\begin{alignat*}{4}
	&R_+(X,Y),\qquad &&S_+(U,V),\qquad &&T_+(X,Y),&\\
	&R_-(U,V),\qquad &&S_-(X,Y),\qquad &&T_-(U,V).&
\end{alignat*}
We are using Lemma~\ref{facilitate} and equations \eqref{list_of_RST} to define the involutions in terms of where they send either $(X,Y)$ or $(U,V)$. By defining the involutions in this way, we ensure that they each descend to the real- and imaginary-symplectic structures on ${\rm Orb}(-1)$ and ${\rm Orb}({\rm i})$ appearing in Proposition~\ref{RST}.

For the next step, we must express each of these 6 involutions in terms of $(Q,P)$, $(X,Y)$, and~$(U,V)$. Once we have these expressions, we will be able to use the proof of Theorem~\ref{HK_red_thm} to track how they descend to involutions $\widetilde{R}_\pm$, $\widetilde{S}_\pm$, $\widetilde{T}_\pm$ on $\tilde{M}$. We claim that
\begin{itemize}\itemsep=0pt
	\item $\widetilde{R}_+$ is identically $\overline{\rho^*}$ on ${\rm Orb}(-1)$ and the imaginary cotangent-lift of the real structure on $\Gr_\C$ with fixed-point set $\Gr_\R$. Conversely, $\widetilde{R}_-$ is identically $-\overline{\rho^*}$ on ${\rm Orb}({\rm i})$ and the real-cotangent lift of the real structure $\Gr_\R\subset\Gr_\C$.
	\item $\widetilde{S}_+$ is equal to $\widetilde{S}_-$ and identically gives $\overline{\sigma}^*$ on ${\rm Orb}({\rm i})$ and $-\overline{\sigma^*}$ on ${\rm Orb}(-1)$.
	\item $\widetilde{T}_+$ is identically $\overline{\tau^*}$ on ${\rm Orb}(-1)$ and the imaginary cotangent-lift of the real structure on $\Gr_\C$ with fixed-point set $\Gr_\Q$. Conversely, $\widetilde{T}_-$ is identically $-\overline{\tau^*}$ on ${\rm Orb}({\rm i})$ and the real-cotangent lift of the real structure $\Gr_\Q\subset\Gr_\C$.
\end{itemize}
It follows that the fixed-point sets for each of these involutions are simultaneously real-symplectic forms with respect to one complex structure, and imaginary-symplectic forms with respect to another. Therefore, from Proposition~\ref{branes} we conclude that they are branes in $\tilde{M}$.

The task of establishing the previous claim is a matter of computation from which we shall spare the reader. We will however flesh out in more detail the resulting complex-Lagrangian submanifolds. Take for instance $R_+$. This involution sends $(U,V)$ to $\bigl(-{\rm i}V^{\mathsf{T}},{\rm i}U^{\mathsf{T}}\bigr)$, and hence, sends $\xi=UV$ in ${\rm Orb}({\rm i})$ to $\xi^{\mathsf{T}}$. Similarly, for $T_+$ one can show that $\xi$ is sent to $-\mathbb{J}\xi^\dagger\mathbb{J}$. The~fixed-point sets are identified with the following submanifolds of ${\rm Orb}(\zeta)\cong\D_\C$.

\begin{Definition}
	The complex symmetric space $\D_\C^\text{ortho}\subset\D_\C$ is the space of decompositions which are orthogonal with respect to the non-degenerate symmetric bilinear form $(x,y)\mapsto x^{\mathsf{T}}y$. Additionally, the set of orthogonal decompositions with respect to the non-degenerate skew-symmetric bilinear form $(x,y)\mapsto x^{\mathsf{T}}\mathcal{J}y$ is denoted by $\D_\C^\text{symp}\subset\D_\C$.
\end{Definition}

Finally, we consider the complex-Lagrangian submanifold of $T^*_{1,0}\Gr_\C$ arising from ${S}_\pm$. This involution sends $(Q,P)$ to $(\Theta Q, -P\Theta)$. Peering back inside the proof of Theorem~\ref{HK_red_thm}, we recall that $\Pi=\Imag Q$ and $\eta=QP$. The fixed-point set is then equal to the space of $\Theta$-invariant $m$-dimensional subspaces $\Pi$ of $\C^n$ with tangent vectors $\eta$ satisfying $\eta=-\Theta\eta\Theta$. Since $\Theta^2=\Id_n$ and $\Theta^\dagger=\Theta$, $\C^n$ decomposes orthogonally into the $\pm1$-eigenspaces $E_{+}\oplus E_{-}$ of $\Theta$. With respect to this decomposition, the fixed-point set is equal to the space $\mathcal{X}$ defined as follows.
\begin{Definition}
	Let $\C^n=E_{+}^{(l)}\oplus E_{-}^{(n-l)}$ be an orthogonal decomposition with respect to the standard Hermitian form $\langle\,,\,\rangle$. Consider the submanifold $\Gr(E_{+})\times\Gr(E_{-})\subset\Gr_\C$ of $m$-dimensional subspaces $\Pi$ which distribute over the decomposition, that is to say, where ${\Pi=\Pi_{+}\oplus\Pi_{-}}$ for $\Pi_{+}\subset E_+$ and $\Pi_-\subset E_-$. The orthogonal complement also distributes over the decomposition $\Pi^\perp=\Pi^\perp_+\oplus\Pi^\perp_-$ and we may identify cotangent vectors to $\Pi$ with linear maps $\eta\colon \Pi^\perp\rightarrow\Pi$. The~space $\mathcal{X}$ is the subbundle of $T^*_{1,0}\Gr_\C$ over $\Gr(E_+)\times\Gr(E_-)$ consisting of those covectors~$\eta$%
\[
		\begin{tikzpicture}
			\node at (0,.5){$\Pi^\perp_+\oplus \Pi^\perp_-$};
			\node at (0,-.5){$\Pi_+\oplus \Pi_-$};
			\draw[-{Latex}] (-0.3,.3) -- ++(-40:.8);
			\draw[-{Latex}] (0.3,.3) -- ++(-140:.8);
		\end{tikzpicture}
\]
	for which $\eta(\Pi^\perp_+)\subset\Pi_{-}$ and $\eta(\Pi^\perp_-)\subset\Pi_{+}$.
\end{Definition}

\begin{Remark}
Observe that $\mathcal{X}$ will typically consist of multiple connected components corresponding to the various possible dimensions of $\Pi^+$. There are two special cases worth noting. When $l=n$, $\mathcal{X}$ is the zero section of $T^*_{1,0}\Gr_\C$. On the other hand, if $l=m$, one of the connected components of $\mathcal{X}$ is the fibre of $T^*_{1,0}\Gr_\C$ over the point $\Pi=E_+$. We recognise that these are indeed complex-Lagrangian submanifolds.
\end{Remark}

\begin{Theorem}The hyperk\"{a}hler reduced space $\tilde{M}$ admits the brane submanifolds listed in the rows of Table~$\ref{RST_table}$. In particular, the complex symmetric spaces $\D_{\C}^{\rm ortho}$ and $\D_\C^{\rm symp}$ are complex-Lagrangian submanifolds of ${\rm Orb}(\zeta)$, and the $\mathcal{X}$ are complex-Lagrangian submanifolds of $T^*_{1,0}\Gr_\C$.
\end{Theorem}

\begin{Corollary}
	The map $\Phi\colon{\rm Orb}({\rm i})\rightarrow T^*_{1,0}\Gr_\C$ restricted to the fixed-point sets of $\widetilde{R}_\pm$, $\widetilde{S}_\pm$, and $\widetilde{T}_\pm$ define real symplectomorphisms:
	\begin{alignat*}{3}
		&(\D_\R,\Imag\Omega_{\rm KKS})\cong(T^*\Gr_\R,\omega_{\rm can}),\qquad && (\D_\C^{\rm ortho},\omega_\textup{K})\cong({\rm i}T^*\Gr_\R,\omega_{\rm can}),& \\
		&(\D_\Q,\Imag\Omega_{\rm KKS})\cong(T^*\Gr_\Q,\omega_{\rm can}),\qquad && (\D_\C^{\rm symp},\omega_\textup{K})\cong({\rm i}T^*\Gr_\Q,\omega_{\rm can}).& \\
		&(\D_\C^{\rm unit},\Real\Omega_{\rm KKS})\cong(\mathcal{X},\omega_\textup{C}),&
	\end{alignat*}
\end{Corollary}

\begin{table}[t]
		\caption{\label{RST_table}Each row represents a brane of $\tilde{M}$. The table follows the template in \eqref{little_table} to demonstrate which submanifolds are real-/imaginary-symplectic forms and complex-Lagrangian submanifolds of $T^*_{1,0}\Gr_\C$ and~${\rm Orb}(\zeta)$.}
$$\arraycolsep=5pt\def\arraystretch{1.6}
		\begin{array}{l|c|c|c}
			&T^*_{1,0}\Gr_\C & {\rm Orb}(-1) & {\rm Orb}({\rm i})\\\hline
			\widetilde{R}_-&T^*\Gr_\R & \D_\C^\text{ortho} & \D_\R \\
			\widetilde{T}_-&T^*\Gr_\Q & \D_\C^\text{symp}& \D_\Q \\ \hline
			\widetilde{S}_\pm & \mathcal{X} & \D_\C^{\rm unit} & \D_\C^{\rm unit} \\ \hline
			\widetilde{R}_+ & {\rm i}T^*\Gr_\R & \D_\R & \D_\C^\text{ortho} \\
			\widetilde{T}_+ & {\rm i}T^*\Gr_\Q & \D_\Q & \D_\C^\text{symp}
		\end{array}
$$
\end{table}

\section{Integrability and dynamical systems}

\subsection{Holomorphic integrability}

\begin{Definition}A collection of holomorphic functions $f_1,\dots,f_n$ on a holomorphic symplectic manifold $M^{(2n)}$ defines a \emph{holomorphic integrable system} if the functions Poisson commute with respect to the complex Poisson bracket and if the derivatives are linearly independent in an open dense subset.
\end{Definition}

The regular fibres of a holomorphic integrable system are complex-Lagrangian submanifolds. These fibres are Lagrangian with respect to both $\omega_R$ and $\omega_I$, and therefore the real and imaginary parts of the $f_k$ define real integrable systems with respect to both $\omega_R$ and $\omega_I$. This is an example of a bi-integrable system and has been observed numerous times before \cite{biham, biham2}.

In the ordinary setting of real integrable systems a great deal of effort can be spent checking a collection of integrals are linearly independent almost everywhere. In this respect working within the rigid category of analytic functions has an advantage: if any two such functions agree in an open set, then they must be identical thanks to the identity theorem. More generally, any subset which enjoys this property is sometimes referred to as a key set \cite{arnoldkeyset,kozlovkeyset}. If we have an analytic real form $N$ of $M$, then every open subset of $N$ is also a key set. This can be seen by working in an analytic chart $\C^{2n}$ in which $N$ appears as $\R^{2n}$ and considering the series expansion. These ideas hold more generally for holomorphic forms, and it is by applying these arguments to $\mathrm{d}f_1\wedge\cdots\wedge \mathrm{d}f_k$ that we prove the following.

\begin{Lemma}\label{keyset_lemma}
	Let $f_1,\dots,f_k$ be holomorphic functions defined on a connected, complex manifold~$M$. If the $\mathrm{d}f_j$ are linearly independent somewhere, then they are independent everywhere in an open dense subset of $M$. If, in addition, $N$ is an analytic real form of $M$, then the $\mathrm{d}f_j$ are also linearly independent in an open dense subset of~$N$.
\end{Lemma}

\begin{Proposition}
Suppose $u_1,\dots,u_n$ is an integrable system of analytic functions on an analytic real-symplectic form $(N,\widehat{\omega}_R)$. The holomorphic extensions of these functions defined in a~neighbourhood of $N$ is a holomorphic integrable system.
\end{Proposition}
\begin{proof}
	Let $f_1,\dots,f_n$ be the holomorphic extensions defined in some neighbourhood of $N$. Since the $f_j$ are purely real on $N$ their Hamiltonian vector fields are tangent to $N$, and so since $N$ is Lagrangian with respect to $\omega_I$, we have
	\[
	\{f_j,f_k\}(p)=\Omega(X_{f_j},X_{f_k})=\omega_R(X_{u_j},X_{u_k})=\{u_j,u_k\}_R(p)=0
	\]
	for $p\in N$. As $N$ is a key set, it follows that $\{f_j,f_k\}$ must be zero everywhere in the neighbourhood of $N$. Finally, if the $\mathrm{d}u_j$ are independent at $p\in N$, then so are the $\mathrm{d}f_j$.
\end{proof}

We now consider a partial converse to this proposition. Suppose $f_1,\dots,f_n$ is a holomorphic integrable system on $M$. By restricting the real and imaginary parts of these functions to a~connected and analytic real-symplectic form~$N$, we obtain a collection of $2n$ real functions. Let~${A}$ denote the algebra of functions on $N$ generated by taking the Poisson bracket with respect to $\widehat{\omega}_R$ between these functions. Fix some $x\in N$ for which the dimension of
\[
\Span\{dw_x\mid w\in{A}\}
\]
is maximal. We can then select $k$ functions $w_1,\dots,w_k$ belonging to ${A}$ whose derivatives form a~basis of this space at $x$. These functions are analytic, and so, since their derivatives are linearly independent at some point they must be independent on an open dense set $U$ of $N$. From the maximality assumption, the fibres of $w_1,\dots,w_k$ must coincide with the fibres of ${A}$ in $U$. The~fibres of ${A}$ are contained in the intersection of the level sets of $f_1,\dots,f_n$ with~$N$. However, as the regular fibres of the holomorphic integrable system are Lagrangian with respect to $\omega_R$, and as $N$ is symplectic with respect to $\omega_R$, any submanifold contained in these intersections must be isotropic with respect to $(N,\widehat{\omega}_R)$. The algebra $A$ is therefore a complete algebra on~$(N,\widehat{\omega}_R)$ in the sense defined in~\cite{bols}.

This construction might appear underwhelming. Indeed, as the regular fibres of $f_1,\dots,f_n$ and $N$ are both real $2n$-dimensional submanifolds in $M$, the generic transversal intersection between them is a point. A foliation of $N$ into points qualifies as a non-commutative integrable system, albeit not an interesting one. On the other hand, it is important to appreciate the significance of this to dynamics. Suppose the flow of a holomorphic Hamiltonian on $M$ admits a holomorphic integrable system $f_1,\dots,f_n$. If the flow leaves $N$ invariant, then ${A}$ is a complete algebra of first integrals. A result of \cite{bols} shows that this implies integrability in the standard sense.
\begin{Theorem}
Suppose the flow of a holomorphic Hamiltonian on $(M,\Omega)$ admits a holomorphic integrable system. If the flow leaves an analytic real-symplectic form $N$ invariant, then the corresponding real Hamiltonian system on $(N,\widehat{\omega}_R)$ is integrable.
\end{Theorem}

\subsection{Compatible momentum}

As we remarked earlier, the level sets of an arbitrary holomorphic integrable system $f_1,\dots, f_n$ and a real-symplectic form will typically intersect transversally in a point. In practice, it is reasonable to expect some compatibility between $f=(f_1,\dots, f_n)\rightarrow\C^n$ and a real-symplectic structure $R$. Ideally, we would like $R$ to act on the fibres of $\mu$. This can be ensured with the following definition.
\begin{Definition}
	Let $\mathfrak{g}$ be a complex Lie algebra and $\mu\colon M\rightarrow\mathfrak{g}^*$ a holomorphic momentum map on the holomorphic symplectic manifold $(M,\Omega)$. The momentum $\mu$ is $R$-\emph{compatible} with respect to a real-symplectic structure $R$ if there exists a real-Poisson structure $\overline{\rho^*}$ on $\mathfrak{g}^*$ satisfying
	\begin{equation}\label{compatiblemomentum}
		{\mu\circ R}=\overline{\rho^*}\circ\mu.
	\end{equation}
\end{Definition}
Notice that this is a generalisation of Definition~\ref{compat_Ham} of what it meant for a Hamiltonian to be $R$-compatible. Now suppose $\mu$ is a $R$-compatible holomorphic integrable system on $M$ and let $\xi_1,\dots, \xi_n$ be a real basis for $\mathfrak{g}^\rho$. The holomorphic integrable system given by $\widehat{\mu}=(g_1,\dots, g_n)$ where $g_k=\langle\mu,\xi_k\rangle$ is functionally equivalent to $\mu$ and purely real on $M^R$. The fixed-point set of a real structure is analytic \cite{analyticfix} and so by Lemma~\ref{keyset_lemma} the integrals $g_1,\dots,g_n$ must be linearly independent in an open dense subset of $M^R$. It follows that the restriction of $\mu$ defines a real integrable system on $M^R$. This generalises the result of \cite{bulgarian1} to holomorphic symplectic manifolds.

\begin{Theorem}\label{real_int_system}
	Let $\mu\colon M\rightarrow\C^n$ be a holomorphic integrable system on a holomorphic symplectic manifold $(M,\Omega)$. If $\mu$ is $R$-compatible with respect to a real-symplectic structure $R$ on $M$ with non-empty fixed-point set $M^R$, then the restriction $\widehat{\mu}\colon M^R\rightarrow\fix\overline{\rho^*}$ is a real integrable system on $\bigl(M^R,\widehat{\omega}_R\bigr)$.
\end{Theorem}

\begin{Proposition}\label{real_momentum}
	Suppose an action of a complex Lie group $G$ on $(M,\Omega)$ by holomorphic symplectomorphisms admits an equivariant momentum map. The action admits an equivariant, $R$-compatible holomorphic momentum map if and only if the $G$-action is $R$-compatible with respect to some real group structure $\rho$. Moreover, when this holds, the restriction $\widehat{\mu}\colon M^R\rightarrow(\mathfrak{g}^\rho)^*$ is an equivariant momentum map for the action of $G^\rho$ on $\bigl(M^R,\widehat{\omega}_R\bigr)$.
\end{Proposition}
\begin{proof}
	If the action of $G$ is $R$-compatible, then differentiating \eqref{equivariant_condition} shows
	\begin{equation}\label{inf_compat}
		R_*(X_\xi)=X_{\rho_*\xi},
	\end{equation}
	where $X_\xi$ denotes the Hamiltonian vector field generated by $\xi\in\mathfrak{g}$. This implies that the~momentum map satisfies
	\[
	{\mu\circ R}=\overline{\rho^*}\circ\mu+c
	\]
	for some constant $c\in\mathfrak{g}$. The equivariance of $\mu$ along with the fact that the coadjoint action is $\rho$-compatible imply that $c$ is central; that is, it is fixed by the coadjoint action. Furthermore, since~$R$ is an involution it follows that $\overline{\rho^*}c=-c$. The map $\mu-c/2$ is then an $R$-compatible equivariant momentum map. Conversely, if \eqref{compatiblemomentum} holds, then so does \eqref{inf_compat}. This can be integrated to show that the $G$-action is $R$-compatible with respect to $\rho$.
\end{proof}

\subsection{Compact real forms of mechanical systems}
	
	A mechanical system on a configuration space $M$ is a Hamiltonian system on the cotangent bundle $T^*M$ with Hamiltonian
	\[
	H(\eta)=\underbrace{g(\eta,\eta)}_\text{kinetic energy}+\underbrace{V(x)}_\text{potential},
	\]
	where $g$ is a fibrewise inner product and $V$ is a function defined on the base $M$. Observe that for any mechanical system the Hamiltonian is invariant under negation $\eta\mapsto-\eta$ of the fibres.
	
	\begin{Theorem}\label{main_thm}
		Let $H^\C$ be a holomorphic function on ${\rm Orb}({\rm i})$ with
		\[
		H^\C\big|_{\D_\mathbb{F}}=H\circ\Phi\big|_{\D_\mathbb{F}},
		\]
		where $H$ is a real Hamiltonian on $T^*\Gr_\mathbb{F}$ and $\D_\mathbb{F}$ is the imaginary-symplectic form for either $\mathbb{F}=\R$ or $\Q$. If $H(\eta)=H(-\eta)$, then $H^\C$ is purely real on the compact real-symplectic form~${\D_{\C}^{\rm unit}\cong\Gr_\C}$.
		
		Secondly, let $H^\C$ be a holomorphic function on ${\rm Orb}({\rm i})\times{\rm Orb}({\rm i})$ with
		\[
		H^\C\big|_{{\rm Orb}({\rm i})}=H\circ\Phi\big|_{{\rm Orb}({\rm i})},
		\]
		where $H$ is a real Hamiltonian on $T^*\Gr_\C$ and ${\rm Orb}({\rm i})$ is an imaginary-symplectic form inside ${\rm Orb}({\rm i})\times{\rm Orb}({\rm i})$. If $H(\eta)=H(-\eta)$, then $H^\C$ is purely real on the compact real-symplectic form~$\Gr_\C\times\Gr_\C$.
	\end{Theorem}
	\begin{proof}
		For when $l=n$, the Hermitian form used in the definition of $\widetilde{S}$ is positive definite and the fixed-point set $\D_{\C}^{\rm unit}$ is identified with $\Gr_\C$. The involution $\widetilde{S}$ is identified via $\Phi$ with the involution which negates covectors $\eta\mapsto -\eta$. As a holomorphic function is determined uniquely by its restriction to a real form, if $H$ is invariant with respect to this involution, then $H^\C\circ\widetilde{S}=\overline{H^\C}$. Hence, $H^\C$ is real on $\fix\widetilde{S}$.
		
		We now consider the twisted involution $(\xi_1,\xi_2)\mapsto (\widetilde{R}_-(\xi_2),\widetilde{R}_-(\xi_1))$ defined on ${\rm Orb}({\rm i})\times{\rm Orb}({\rm i})$ with the holomorphic symplectic form $\frac{1}{2}\Omega_{\rm KKS}\oplus\frac{1}{2}\Omega_{\rm KKS}$. The fixed-point set
		\begin{equation}\label{conj_diag}
			\big\{\bigl(\xi,\widetilde{R}_-(\xi)\bigr)\mid \xi\in{\rm Orb}({\rm i})\big\}
		\end{equation}
		is an imaginary-symplectic form symplectomorphic to $({\rm Orb}({\rm i}),\Imag\Omega_{\rm KKS})$. The involutions $\widetilde{R}_-$ and $\widetilde{S}$ commute, and hence, $H(\eta)=H(-\eta)$ implies $H^\C\circ\bigl(\widetilde{S}\oplus\widetilde{S}\bigr)=\overline{H^\C}$. It follows that $H^\C$ is real on $\fix\bigl(\widetilde{S}\oplus\widetilde{S}\bigr)$.
	\end{proof}

\begin{Remark}The theorem above is a more precise version of our main theorem appearing in the introduction. To apply the theorem, we suppose we have a mechanical system on $T^*\Gr_\mathbb{F}$ given by a real analytic Hamiltonian $H$. The complexification $H^\C$ is a holomorphic function which is not necessarily defined everywhere. For this reason, there can be no guarantee that the restriction to the compact real form is globally defined. It must however be defined on some non-empty open subset since the intersection of the imaginary-symplectic form $\D_\mathbb{F}$ and the compact real-symplectic form is non-empty. In the examples to follow, we shall see that the Hamiltonians are defined everywhere on the compact form.
\end{Remark}
	We will now demonstrate a few applications of this result. Owing to the low-dimensional isomorphisms $\RP^1\cong S^1$, $\CP^1\cong S^2$, and $\RP^3\cong {\rm SO}(3)$, we can use the theorem to find compact real forms of well-known mechanical systems: the simple pendulum, the spherical pendulum, and the rigid body. To do this, we will need to see how kinetic energy and potential energy pull back through $\Phi$.
	\begin{Proposition}\label{pullbacks_explicit}
		The standard ${\rm U}(n)$-invariant kinetic energy $|\eta|^2$ on $T^*\CP^{n-1}$ pulls back through $\Phi$ to the function $|\xi|^2-1$ on ${\rm Orb}({\rm i})$, where $|\xi|^2=\Tr\xi\xi^\dagger$.
		
		Let $\pi\colon T^*\CP^{n-1}\rightarrow\CP^{n-1}$ denote projection onto the base. The map $\widetilde{\pi}\colon{\rm Orb}({\rm i}) \rightarrow\fix\widetilde{S}=\D_{\C}^{\rm unit}$ given by
		\begin{equation}\label{bundle_map}
			\widetilde{\pi}(\xi)=\frac{1}{2(|\xi|+1)}\biggl[\frac{{\rm i}\bigl(\xi\xi^\dagger+\xi^\dagger\xi\bigr)}{|\xi|}+\xi-\xi^\dagger\biggr]
		\end{equation}
		satisfies $\Phi\circ\widetilde{\pi}=\pi\circ\Phi$.
	\end{Proposition}
	\begin{proof}
		A covector $\eta\colon\C^n/\Pi\rightarrow\Pi$ in the standard metric has magnitude $\Tr\eta\eta^\dagger$. Using the expressions in Proposition~\ref{Phi_for_CPn}, we obtain the first part directly by calculation.
		
		For the second part, $\widetilde{\pi}(\xi)$ has image $\Span\{z\}$ and kernel $(\Span\{z\})^\perp$ for $z$ in \eqref{z_w_defn}. Using~\eqref{xi_as_xy}, we see that $\widetilde{\pi}(\xi)$ must be equal to ${\rm i}{z}{z}^\dagger/|z|^2$. Another calculation using Proposition~\ref{Phi_for_CPn} gives the expression above.
	\end{proof}

\begin{Example}[the simple pendulum]
		For $m=1$ and $n=2$, the orbit ${\rm Orb}({\rm i})$ may be written~as%
		\begin{equation}\label{2by2}
			\biggl\{\xi=\frac{1}{2}\begin{pmatrix}
				t+{\rm i}x & y+{\rm i}z \\ -y+{\rm i}z & t-{\rm i}x
			\end{pmatrix} \mid t={\rm i},\,\det\xi=0\biggr\}\subset\mathfrak{gl}_2\C^*.
		\end{equation}
		This gives the affine variety $x^2+y^2+z^2=1$ over $\C^3$ which we denote by $\CS^2$ and call the complex 2-sphere. Writing $\text{x}=(x,y,z)$, we find that
		\begin{equation}\label{R_twiddle}
\widetilde{R}_-(\text{x})=(\overline{x},-\overline{y},\overline{z}),\qquad\text{and}\qquad\widetilde{S}(\text{x})=(\overline{x},\overline{y},\overline{z}).
		\end{equation}
		Furthermore, the map $\widetilde{\pi}$ in \eqref{bundle_map} simplifies considerably to
		\begin{equation}\label{simple_pi}
			\widetilde{\pi}(\text{x})=\frac{\sqrt{2}}{\sqrt{1+|\text{x}|^2}}\Real(\text{x}),
		\end{equation}
		where $|\text{x}|^2=|x|^2+|y|^2+|z|^2$. We claim that the holomorphic Hamiltonian
		\begin{equation*}
			H^\C(\text{x})=\frac{1}{2}\bigl(x^2-y^2+z^2-1\bigr)+\frac{x\sqrt{2}}{\sqrt{1+x^2-y^2+z^2}}
		\end{equation*}
		restricted to $\fix\widetilde{R}_-=\D_\R$ is the pullback of the usual Hamiltonian for the simple pendulum through the symplectomorphism $(\D_\R,\Imag\Omega_{\rm KKS})\overset{\Phi}{\rightarrow}\bigl(T^*\RP^1,\omega_\text{can}\bigr)$. The kinetic energy comes from Proposition~\ref{pullbacks_explicit}, where $|\eta|^2=|\xi|^2-1=\bigl(|\text{x}|^2-1\bigr)/2$. For the potential, we note that the base $\RP^1$ is identified with the circle $x^2+z^2=1$, and so we take the potential to be the $x$-component in \eqref{simple_pi}. Restricted to the real-symplectic form $\fix\widetilde{S}=S^2$ this is
		\begin{equation*}
			H^\C\big|_{S^2}=-\cos^2\psi+\cos\varphi
		\end{equation*}
		given in terms of spherical polar coordinates. Compare this with the original Hamiltonian $H(
		\eta,\theta)=|\eta|^2+\cos\theta$ on $T^*S^1$.
	\end{Example}
	
\begin{figure}[t]
\centering\includegraphics[scale=0.8]{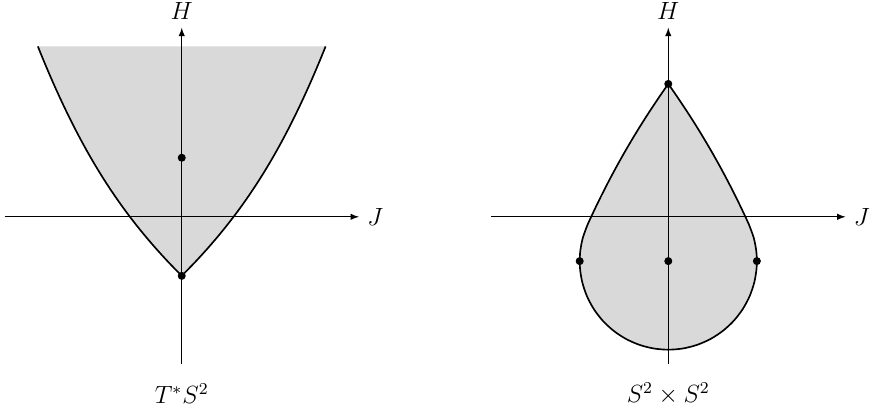}
\caption{\label{diagram}Energy-momentum diagram for the spherical pendulum and its `compact real form'. The black points are rank zero critical points of $(H,J)$.}
	\end{figure}

	\begin{Example}[the spherical pendulum]
		Consider the pullback of the Hamiltonian for the spherical pendulum through the symplectomorphism \smash{$\bigl(\CS^2,\Imag\Omega_{\rm KKS}\bigr)\overset{\Phi}{\rightarrow}\bigl(T^*\CP^1,\omega_\text{can}\bigr)$}. As in the previous example, the kinetic energy pulls back to $\bigl(|\text{x}|^2-1\bigr)/2$. If we take the vertical to be about $(0,1,0)$, then the $y$-component in \eqref{simple_pi} gives the height of the pendulum. Thus,
		\[
		H\circ\Phi\big|_{\CS^2}=\frac{1}{2}\bigl(|x|^2+|y|^2+|z|^2-1\bigr)+\frac{y+\overline{y}}{\sqrt{2}\sqrt{1+|x|^2+|y|^2+|z|^2}}.
		\]
		The Hamiltonian is invariant under the ${\rm U}(1)$-action on $T^*\CP^1$ which rotates about the vertical. Since $\Phi$ is ${\rm U}(2)$-equivariant, this corresponds to the circle action on $\CS^2$ which fixes $x$ and $z$. As $\CS^2$ is equipped with the imaginary part of $\Omega_{\rm KKS}$, this circle action is the Hamiltonian flow generated by the imaginary part of $y$. Hence, the angular momentum $J$ about the $y$-axis pulls pack to
		\[
		J\circ\Phi=\frac{y-\overline{y}}{2{\rm i}}.
		\]
		If we write elements in $\CS^2\times\CS^2$ as pairs $(\text{x}_1,\text{x}_2)$, then using \eqref{R_twiddle} and the definition of the~conjugate-diagonal copy of $\CS^2$ in \eqref{conj_diag} we extend $H$ and $J$ to holomorphic functions
\begin{align*}
&H^\C=\frac{1}{2}(x_1x_2-y_1y_2+z_1z_2-1)+\frac{y_1-y_2}{\sqrt{(x_1+x_2)^2+(y_1-y_2)^2+(z_1+z_2)^2}},\\
&J^\C=\frac{(y_1+y_2)}{2{\rm i}}.
\end{align*}
This integrable system is compatible with respect to the real structure $(H,J)\mapsto (\overline{H},-\overline{J})$ on~$\C^2$. Therefore, from Theorem~\ref{real_int_system} it provides a real integrable system on the compact~real form~${S^2\times S^2}$. For reference, we provide in Figure~\ref{diagram} the energy-momentum diagrams for the~systems on $T^*S^2$ and $S^2\times S^2$.
\end{Example}
	
\begin{Example}[the rigid body]The matrix in \eqref{2by2} allows us to identify real vectors $(t,x,y,z)$ in $\R^4\setminus\{0\}$ with a group of $2\times 2$ matrices. We can use this to define a group structure on $\RP^3$ which turns out to be isomorphic to ${\rm SO}(3)$. The Hamiltonian $H$ for the rigid body is a~left-invariant function on $T^*{\rm SO}(3)$. Consequently, it may be written as $H=f\circ\mu$ where $\mu$ is the momentum map for the right-action of ${\rm SO}(3)$ on $T^*{\rm SO}(3)$ and $f$ is a positive definite quadratic on the Lie algebra $\mathfrak{so}(3)^*\cong\R^3$ which defines the kinetic energy.
		
Now consider $\D_\C$ for the case $m=1$ and $n=4$. As before, we may identify vectors $\text{x}$ in~$\C^4$ with matrices $X$ in $M_2(\C)$ and consider the ${\rm SL}_2\C$-action given by right multiplication. This~gives a Hamiltonian action of ${\rm SL}_2\C\subset {\rm GL}_4\C$ on ${\rm Orb}({\rm i})$ with momentum map,
		\[
		\mu^\C\colon \ {\rm Orb}({\rm i})\longrightarrow\mathfrak{sl}_2\C^*.
		\]
		This action is both $\widetilde{R}_-$ and $\widetilde{S}$-compatible for the real group structure $\sigma(g)=g^{-\dagger}$ with $\fix\sigma={\rm SU}(2)$. It follows from the ${\rm U}(4)$-equivariance of $\Phi$ that the $\fix\sigma={\rm SU}(2)$-action on $\D_\R$ pushes forward to give the action of ${\rm SO}(3)$ on $T^*{\rm SO}(3)$ mentioned above. Therefore, from Proposition~\ref{real_momentum}
		\[
		\mu^\C\big|_{\D_\R}=\mu\circ\Phi\big|_{\D_\R}.
		\]
		The holomorphic momentum $\mu^\C$ is also compatible with respect to $\widetilde{S}$ and so restricts to the momentum map for the ${\rm SU}(2)$-action on $\D_{\C}^{\rm unit}\cong\CP^3$ which is given by
		\begin{equation*}
			\mu^\C\big|_{\CP^3}\colon \ \CP^3\longrightarrow\mathfrak{su}(2)^*,\qquad [\text{x}]\longmapsto {\rm i}\biggl(\frac{X^\dagger X}{|\text{x}|^2}-\Id\biggr).
		\end{equation*}
		Thus, we obtain a holomorphic Hamiltonian $f\circ\mu^\C$ on ${\rm Orb}({\rm i})$ which is purely real on the real-symplectic form $\CP^3$ and restricts to the left-invariant rigid-body Hamiltonian on the imaginary-symplectic form \smash{$(\D_\R,\Imag\Omega_{\rm KKS})\overset{\Phi}{\rightarrow}(T^*\RP^3,\omega_\text{can})$}.
	\end{Example}

	\section{Concluding comments and scope for further work}
	Our work in this paper has relied upon the existence of a hyperk\"{a}hler structure on a particular class of coadjoint orbits. More generally, for any coadjoint orbit $G/H$ of a compact Lie group~$G$, the complexified orbit $G^\C/H^\C$ admits a hyperk\"{a}hler structure \cite{french2,kovalev}. Moreover, such orbits are known to be diffeomorphic to $T^*(G/H)$. Therefore, it seems reasonable to suspect that a~more general version of Theorem~\ref{main_thm} might hold for mechanical systems on any coadjoint orbit of a~compact Lie group.
	
Another example of hyperk\"{a}hler geometry appearing in classical mechanics is the Calogero--Moser system for $n$ mutually interacting point particles on a line. Here it is known that the collisions can be regularised to obtain a symplectic reduced space for the action of ${\rm U}(n)$ on~$T^*\mathfrak{u}(n)$~\cite{calogeromoser}. The complexified Calogero--Moser space is also hyperk\"{a}hler, and by changing the complex structure the resulting holomorphic symplectic manifold is the Hilbert scheme of~$n$ points in the affine plane \cite{george}. As with our examples above, one could try to find additional real forms for the Calogero--Moser system by considering the real forms of ${\rm GL}_n\C$ or by looking for complex-Lagrangian submanifolds of the Hilbert scheme. In this regard, the literature concerning involutions of holomorphic symplectic manifolds could be of use \cite{beau,hilbert_involution, brane_involutions}.
	
Theorem~\ref{real_int_system} suggests that it might be of interest to complexify a real integrable system and look for compatible real forms. The existence of a compact integrable real form is of particular interest since their theory is better understood. In efforts towards the classification of integrable systems the nature of the critical points of the momentum map plays a significant role. Non-degenerate critical points are determined locally by a Cartan subalgebra of the real symplectic Lie algebra \cite{int_theory}. In this respect, the study of holomorphic integrable systems might afford an advantage over its real counterpart: unlike the real symplectic Lie algebra, every Cartan subalgebra of the complex symplectic Lie algebra is equivalent up to conjugacy. This implies that for a holomorphic integrable system, all non-degenerate rank zero critical points are locally the same. It is intriguing to contemplate how this might contribute to the task of classification, similar in spirit to that in \cite{anton}.
	
Finally, our original motivation for considering real forms was to translate the study of one dynamical system into another via the complexification; somewhat analogous to the Wick rotation used in physics. In this vein, we would like to report on the first author's recent work at transferring the study of the 2-body problem on a sphere to the problem on the hyperboloid \cite{paper3plus}. For this particular system, the classification of relative equilibria on the sphere can be used to obtain a classification for those on the hyperboloid by working on the complexified holomorphic system. This prompts us to ask more generally how knowledge of one real form might be used to deduce properties for another.
	
\subsection*{Acknowledgements}
At the time of writing, Philip Arathoon was funded by an EPSRC Doctoral Prize Award hosted by the University of Manchester and Marine Fontaine was supported by the FWO-EoS Project G0H4518N.
We would like to extend special thanks to James Montaldi for carefully reading an early draft and for suggesting many helpful comments throughout its preparation. Additional thanks are owed to the referees whose insights and recommendations helped to improve this paper.\looseness=-1

\pdfbookmark[1]{References}{ref}
\LastPageEnding


\begin{thebibliography}{99}
\footnotesize\itemsep=0pt

\bibitem{Moerbeke}
Adler M., van Moerbeke P., Vanhaecke P., Algebraic integrability, {P}ainlev\'e
 geometry and {L}ie algebras, \textit{Ergeb. Math. Grenzgeb.~(3)}, Vol.~47,
 \href{https://doi.org/10.1007/978-3-662-05650-9}{Springer}, Berlin, 2004.

\bibitem{paper3plus}
Arathoon P., Unifying the hyperbolic and spherical {$2$}-body problem with
 biquaternions, \href{https://doi.org/10.1134/S1560354723060011}{\textit{Regul. Chaotic Dyn.}} \textbf{28} (2023), 822--834,
 \href{https://arxiv.org/abs/2012.12166}{arXiv:2012.12166}.

\bibitem{arnoldkeyset}
Arnold V.I., Kozlov V.V., Neishtadt A.I., Dynamical systems.~III. Mathematical
 aspects of classical and celestial mechanics, 3rd~ed., \textit{Encyclopaedia Math.
 Sci.}, Vol.~3, \href{https://doi.org/10.1007/978-3-540-48926-9}{Springer}, Berlin, 2006.

\bibitem{biham}
Bates L., Cushman R., Complete integrability beyond {L}iouville--{A}rnol'd,
 \href{https://doi.org/10.1016/S0034-4877(05)80042-4}{\textit{Rep. Math. Phys.}} \textbf{56} (2005), 77--91.

\bibitem{beau}
Beauville A., Antisymplectic involutions of holomorphic symplectic manifolds,
 \href{https://doi.org/10.1112/jtopol/jtr002}{\textit{J.~Topol.}} \textbf{4} (2011), 300--304, \href{https://arxiv.org/abs/1008.3108}{arXiv:1008.3108}.

\bibitem{french2}
Biquard O., Sur les \'equations de {N}ahm et la structure de {P}oisson des
 alg\`ebres de {L}ie semi-simples complexes, \href{https://doi.org/10.1007/BF01446293}{\textit{Math. Ann.}} \textbf{304}
 (1996), 253--276.

\bibitem{french}
Biquard O., Gauduchon P., G\'eom\'etrie hyperk\"ahl\'erienne des espaces
 hermitiens sym\'etriques complexifi\'es, Universit\'e de Grenoble I, \href{https://doi.org/10.5802/tsg.199}{Institut
 Fourier}, Saint-Martin-d'H\`eres, 1998, 127--173.

\bibitem{bols}
Bolsinov A.V., Jovanovi\'c B., Noncommutative integrability, moment map and
 geodesic flows, \href{https://doi.org/10.1023/A:1023023300665}{\textit{Ann. Global Anal. Geom.}} \textbf{23} (2003),
 305--322, \href{https://arxiv.org/abs/math-ph/0109031}{arXiv:math-ph/0109031}.

\bibitem{calabi_metric}
Calabi E., M\'etriques k\"ahl\'eriennes et fibr\'es holomorphes, \href{https://doi.org/10.24033/asens.1367}{\textit{Ann.
 Sci. \'{E}cole Norm. Sup.}} \textbf{12} (1979), 269--294.

\bibitem{hilbert_involution}
Cattaneo A., Automorphisms of {H}ilbert schemes of points on a generic
 projective {K}3 surface, \href{https://doi.org/10.1002/mana.201800557}{\textit{Math. Nachr.}} \textbf{292} (2019),
 2137--2152, \href{https://arxiv.org/abs/1801.05682}{arXiv:1801.05682}.

\bibitem{crainic}
Crainic M., Fernandes R.L., Integrability of {P}oisson brackets,
 \href{https://doi.org/10.4310/jdg/1090415030}{\textit{J.~Differential Geom.}} \textbf{66} (2004), 71--137,
 \href{https://arxiv.org/abs/math.DG/0210152}{arXiv:math.DG/0210152}.

\bibitem{pcrook}
Crooks P., Rayan S., Abstract integrable systems on hyperk\"ahler manifolds
 arising from {S}lodowy slices, \href{https://doi.org/10.4310/MRL.2019.v26.n1.a2}{\textit{Math. Res. Lett.}} \textbf{26} (2019),
 9--33, \href{https://arxiv.org/abs/1706.05819}{arXiv:1706.05819}.

\bibitem{biham2}
Doss-Bachelet C., Fran\c{c}oise J.P., Integrable {H}amiltonian systems with two
 degrees of freedom associated with holomorphic functions, \href{https://doi.org/10.1007/BF02551194}{\textit{Theoret.
 and Math. Phys.}} \textbf{122} (2000), 170--175.

\bibitem{brane_involutions}
Franco E., Jardim M., Menet G., Brane involutions on irreducible holomorphic
 symplectic manifolds, \href{https://doi.org/10.1215/21562261-2018-0009}{\textit{Kyoto~J. Math.}} \textbf{59} (2019), 195--235,
 \href{https://arxiv.org/abs/1606.09040}{arXiv:1606.09040}.

\bibitem{witten2}
Gaiotto D., Witten E., Probing quantization via branes, in Surveys in
 Differential Geometry 2019. {D}ifferential Geometry, {C}alabi--{Y}au theory,
 and General Relativity. {P}art 2, \textit{Surv. Differ. Geom.}, Vol.~24,
 \href{https://doi.org/10.4310/SDG.2019.v24.n1.a8}{International Press}, Boston, MA, 2022, 293--402, \href{https://arxiv.org/abs/2107.12251}{arXiv:2107.12251}.

\bibitem{bulgarian3}
Gerdjikov V.S., Kyuldjiev A., Marmo G., Vilasi G., Complexifications and real
 forms of {H}amiltonian structures, \href{https://doi.org/10.1140/epjb/e2002-00281-y}{\textit{Eur. Phys.~J.~B Condens. Matter
 Phys.}} \textbf{29} (2002), 177--181.

\bibitem{bulgarian4}
Gerdjikov V.S., Kyuldjiev A., Marmo G., Vilasi G., Construction of real forms
 of complexified {H}amiltonian dynamical systems, in Nonlinear Physics: Theory
 and Experiment,~{II} ({G}allipoli, 2002), \href{https://doi.org/10.1142/9789812704467_0024}{World Scientific}, River Edge, NJ,
 2003, 172--178.

\bibitem{bulgarian1}
Gerdjikov V.S., Kyuldjiev A., Marmo G., Vilasi G., Real {H}amiltonian forms of
 {H}amiltonian systems, \href{https://doi.org/10.1140/epjb/e2004-00158-1}{\textit{Eur. Phys.~J.~B Condens. Matter Phys.}}
 \textbf{38} (2004), 635--649, \href{https://arxiv.org/abs/nlin.SI/0310005}{arXiv:nlin.SI/0310005}.

\bibitem{witten1}
Gukov S., Witten E., Branes and quantization, \href{https://doi.org/10.4310/atmp.2009.v13.n5.a5}{\textit{Adv. Theor. Math. Phys.}}
 \textbf{13} (2009), 1445--1518, \href{https://arxiv.org/abs/0809.0305}{arXiv:0809.0305}.

\bibitem{anton}
Izosimov A., Singularities of integrable systems and algebraic curves,
 \href{https://doi.org/10.1093/imrn/rnw168}{\textit{Int. Math. Res. Not.}} \textbf{2017} (2017), 5475--5524,
 \href{https://arxiv.org/abs/1509.08996}{arXiv:1509.08996}.

\bibitem{calogeromoser}
Kazhdan D., Kostant B., Sternberg S., Hamiltonian group actions and dynamical
 systems of {C}alogero type, \href{https://doi.org/10.1002/cpa.3160310405}{\textit{Comm. Pure Appl. Math.}} \textbf{31}
 (1978), 481--507.

\bibitem{kovalev}
Kovalev A.G., Nahm's equations and complex adjoint orbits,
 \href{https://doi.org/10.1093/qmath/47.1.41}{\textit{Quart.~J.~Math. Oxford Ser.~(2)}} \textbf{47} (1996), 41--58.

\bibitem{kozlovkeyset}
Kozlov V.V., Integrability and nonintegrability in {H}amiltonian mechanics,
 \href{https://doi.org/10.1070/RM1983v038n01ABEH003330}{\textit{Russ. Math. Surv.}} \textbf{38} (1983), 3--67.

\bibitem{kulkarni}
Kulkarni R.S., On complexifications of differentiable manifolds,
 \href{https://doi.org/10.1007/BF01389901}{\textit{Invent. Math.}} \textbf{44} (1978), 46--64.

\bibitem{bulgarian2}
Kyuldjiev A., Gerdjikov V., Marmo G., Vilasi G., Real forms of complexified
 {H}amiltonian dynamics, in Geometry, Integrability and Quantization ({V}arna,
 2001), \textit{Geom. Integrability Quantization}, Vol.~3, \href{https://doi.org/10.7546/giq-3-2002-318-327}{Coral Press
 Scientific Publishing}, Sofia, 2002, 318--327.

\bibitem{int_theory}
Pelayo A., V\~u Ng\d{o}c S., Symplectic theory of completely integrable
 {H}amiltonian systems, \href{https://doi.org/10.1090/S0273-0979-2011-01338-6}{\textit{Bull. Amer. Math. Soc. (N.S.)}} \textbf{48}
 (2011), 409--455, \href{https://arxiv.org/abs/1306.0115}{arXiv:1306.0115}.

\bibitem{analyticfix}
Sepp\"al\"a M., Quotients of complex manifolds and moduli spaces of {K}lein
 surfaces, \href{https://doi.org/10.5186/aasfm.1981.0622}{\textit{Ann. Acad. Sci. Fenn. Ser.~A~I Math.}} \textbf{6} (1981),
 113--124.

\bibitem{realform2}
Shutrick H.B., Complex extensions, \href{https://doi.org/10.1093/qmath/9.1.189}{\textit{Quart.~J.~Math. Oxford Ser.~(2)}}
 \textbf{9} (1958), 189--201.

\bibitem{dualpairs}
Skerritt P., Vizman C., Dual pairs for matrix groups, \href{https://doi.org/10.3934/jgm.2019014}{\textit{J.~Geom. Mech.}}
 \textbf{11} (2019), 255--275, \href{https://arxiv.org/abs/1805.01519}{arXiv:1805.01519}.

\bibitem{realm}
Vanhaecke P., Integrable systems in the realm of algebraic geometry, 2nd~ed.,
 \textit{Lecture Notes in Math.}, Vol.~1638, \href{https://doi.org/10.1007/3-540-44576-5}{Springer-Verlag}, Berlin,
 2001.

\bibitem{realform1}
Whitney H., Bruhat F., Quelques propri\'et\'es fondamentales des ensembles
 analytiques-r\'eels, \href{https://doi.org/10.1007/BF02565913}{\textit{Comment. Math. Helv.}} \textbf{33} (1959),
 132--160.

\bibitem{george}
Wilson G., Collisions of {C}alogero--{M}oser particles and an adelic
 {G}rassmannian (with an appendix by {I.G.}~{M}acdonald), \href{https://doi.org/10.1007/s002220050237}{\textit{Invent.
 Math.}} \textbf{133} (1998), 1--41.

\end{thebibliography}
\end{document}